\documentclass[]{arxiv-tc}

\usepackage[utf8]{inputenc}

\usepackage{amsmath}
\usepackage{amsthm}
\usepackage{amssymb}
\usepackage{graphicx}
\graphicspath{{tikzfigs/}{figs/}}
\usepackage{mathtools}
\usepackage{subcaption}
\usepackage{array}
\usepackage{color}
\usepackage{hhline}
\usepackage{multirow}
\usepackage{enumitem}
\usepackage{multicol}
\usepackage{bm}
\usepackage{ulem}
\usepackage{booktabs}
\usepackage{algorithm}
\usepackage{algpseudocode}
\usepackage{lipsum}
\usepackage{bibentry}

\usepackage[version=4]{mhchem}
\usepackage{siunitx}
\usepackage{longtable,tabularx}
\makeatletter
\newcommand*\bigcdot{\mathpalette\bigcdot@{.4}}
\newcommand*\bigcdot@[2]{\mathbin{\vcenter{\hbox{\scalebox{#2}{$\m@th#1\bullet$}}}}}
\makeatother

\newcommand{\sref}[1]{\S\ref{#1}}
\newcommand{\oursec}[1]{\vspace{0.1cm}\noindent{\normalsize\textit{#1}:}\,}

\newcommand{\boxing}[6]{
	\begin{figure*}[#1]
		\noindent\makebox[\textwidth][c] {
			\fbox{
				\begin{minipage}{#4}
					\vspace{0.2cm}
					\begin{#2} {\bf \textit{#5}} \label{#3} \end{#2}
					#6
				\end{minipage}
			}
		}
	\end{figure*}
}

\newcommand{\bvec}[1]{#1}

\newcommand{\definedas}{\coloneqq}
\newcommand{\dg}{^{\circ}}

\newcommand{\OMEGA}[1]{\Omega\big(#1\big)}
\newcommand{\diag}[1]{\mathop{\mathrm{diag}}\left(#1\right)}

\newcommand{\dotprod}[2]{#1 \bigcdot #2}
\newcommand{\crossprod}[2]{#1 \times #2}
\newcommand{\ith}[2]{#1^{\textit{#2}}}

\newcommand{\intee}[2]{[#1\,,#2]}
\newcommand{\intei}[2]{[#1\,,#2)}
\newcommand{\intie}[2]{(#1\,,#2]}
\newcommand{\intii}[2]{(#1\,,#2)}

\newcommand{\real}{\mathbb{R}}

\newcommand{\spd}[1]{\mathbb{S}^{#1}_{++}}

\newcommand{\twonorm}[1]{\left\|#1\right\|_2}

\newcommand{\refframe}{\mathcal{F}}
\newcommand{\inertial}{\mathcal{I}}
\newcommand{\body}{\mathcal{B}}
\newcommand{\iframe}{\refframe_\inertial}
\newcommand{\bframe}{\refframe_\body}
\newcommand{\ex}{\bvec{e}_{1}}

\newcommand{\qidentity}{\bvec{q}_{\textit{id}}}

\newcommand{\subin}{\textit{in}}
\newcommand{\subig}{\textit{ig}}
\newcommand{\subf}{\textit{f}}
\newcommand{\tin}{t_{\subin}} 
\newcommand{\tig}{t_{\subig}} 
\newcommand{\tf}{t_{\subf}} 
\newcommand{\tc}{t_{\textit{c}}} 
\newcommand{\tb}{t_{\textit{b}}} 
\newcommand{\tcmax}{t_{\textit{c},\textit{max}}} 

\newcommand{\m}{m}

\newcommand{\rI}{\bvec{r}_\inertial}
\newcommand{\vI}{\bvec{v}_\inertial}
\newcommand{\vB}{\bvec{v}_\body}
\newcommand{\qIB}{\bvec{q}_{\body\leftarrow\inertial}}
\newcommand{\omegaB}{\bvec{\omega}_\body}

\newcommand{\TB}{\bvec{T}_\body}

\newcommand{\aeroB}{\bvec{A}_\body}
\newcommand{\xx}{\bvec{x}} 
\newcommand{\uu}{\bvec{u}} 

\newcommand{\rIdot}{\dot{\bvec{r}}_\inertial}
\newcommand{\vIdot}{\dot{\bvec{v}}_\inertial}
\newcommand{\qIBdot}{\dot{\bvec{q}}_{\body\leftarrow\inertial}}
\newcommand{\omegaBdot}{\dot{\bvec{\omega}}_\body}

\newcommand{\rIin}{\bvec{r}_{\inertial,\subin}}
\newcommand{\vIin}{\bvec{v}_{\inertial,\subin}}

\newcommand{\mig}{\m_{\subig}}

\newcommand{\omegaBig}{\bvec{0}}

\newcommand{\prig}{p_{\bvec{r},\subig}}
\newcommand{\pvig}{p_{\bvec{v},\subig}}

\newcommand{\rIf}{\bvec{0}}

\newcommand{\vIf}{\bvec{0}}
\newcommand{\qIBf}{\qidentity}
\newcommand{\omegaBf}{\bvec{0}}

\newcommand{\mdotalpha}{\alpha_{\dot{\m}}}
\newcommand{\mdotbeta}{\beta_{\dot{\m}}}
\newcommand{\rTB}{\bvec{r}_{T,\body}}
\newcommand{\rCPB}{\bvec{r}_{\textit{cp},\body}}
\newcommand{\inertia}{J_\body}
\newcommand{\Ca}{C_{A}}

\newcommand{\Sa}{S_{A}}

\newcommand{\tiltmax}{\theta_{\textit{max}}}
\newcommand{\omegamax}{\omega_{\textit{max}}}
\newcommand{\glideslope}{\gamma_{\textit{gs}}}
\newcommand{\mdry}{m_{\textit{dry}}}
\newcommand{\gimbalmax}{\delta_{\textit{max}}}
\newcommand{\Tmin}{T_{\textit{min}}}
\newcommand{\Tmax}{T_{\textit{max}}}

\newcommand{\stcnx}{{n_{z}}}
\newcommand{\stcng}{{n_{g}}}
\newcommand{\stcnf}{{n_{c}}}
\newcommand{\stcg}{g} 
\newcommand{\stcf}{c} 
\newcommand{\stch}{h} 
\newcommand{\stcs}{\sigma} 
\newcommand{\stcshat}{\hat{\stcs}} 
\newcommand{\stcx}{\bvec{z}} 

\newcommand{\aoa}{\alpha}
\newcommand{\Vaoa}{V_{\aoa}}
\newcommand{\aoamax}{\aoa_{\textit{max}}}

\newcommand{\gaoa}[1]{\stcg_{\aoa}\big(#1\big)}
\newcommand{\faoa}[2]{\stcf_{\aoa}\big(#1,#2\big)}
\newcommand{\haoa}[2]{\stch_{\aoa}\big(#1,#2\big)}
\newcommand{\gor}[2]{\stcg_{\lor,#1}\big(#2\big)}
\newcommand{\for}[1]{\stcf_{\lor}\big(#1\big)}
\newcommand{\hor}[1]{\stch_{\lor}\big(#1\big)}
\newcommand{\gand}[2]{\stcg_{\land,#1}\big(#2\big)}
\newcommand{\fand}[1]{\stcf_{\land}\big(#1\big)}
\newcommand{\hand}[1]{\stch_{\land}\big(#1\big)}

\newcommand{\gI}{\bvec{g}_\inertial}
\newcommand{\density}{\rho}

\newcommand{\cIB}{C_{\body\leftarrow\inertial}}
\newcommand{\cBI}{C_{\inertial\leftarrow\body}}

\newcommand{\Hgs}{H_{\gamma}}
\newcommand{\Htilt}{H_{\theta}}

\newcommand{\nz}{{n_z}} 
\newcommand{\dilation}{s}
\newcommand{\zz}{\bvec{z}} 
\newcommand{\sso}{\bar{\dilation}} 

\newcommand{\KK}{K}

\newcommand{\wvc}{w_{\nu}}
\newcommand{\Wtr}{W_{\textit{tr}}}

\newcommand{\epsvc}{\epsilon_{\textit{vc}}}
\newcommand{\epstr}{\epsilon_{\textit{tr}}}

\newcommand{\LU}{U_L}
\newcommand{\TU}{U_T}
\newcommand{\MU}{U_M}

\usepackage{rotating}
\usepackage{tikz}
\usetikzlibrary{patterns}
\usetikzlibrary{math}
\usetikzlibrary{scopes}
\usetikzlibrary{fadings}
\usetikzlibrary{arrows,calc,decorations.pathmorphing}
\usetikzlibrary{matrix,positioning}
\tikzset{>=latex}

\definecolor{beige}{RGB}{245,245,220}
\definecolor{darkred}{rgb}{0.90,0.00,0.00}
\definecolor{darkgreen}{rgb}{0.00,0.45,0.00}

\tikzset{ 
table/.style={
  matrix of math nodes,
  row sep=-\pgflinewidth,
  column sep=-\pgflinewidth,
  nodes={rectangle,draw=white,text width=6em,align=left},
  text depth=0.25ex,
  text height=2ex,
  nodes in empty cells
  },
  title/.style={font=\large}
}

\definecolor{ucol}{RGB}{255,0,0}
\definecolor{gcol}{RGB}{0,120,0}
\definecolor{scol}{RGB}{63, 226, 45}
\definecolor{vcol}{RGB}{0,0,0}
\definecolor{acol}{RGB}{0, 128, 255}
\definecolor{dcol}{RGB}{204,102,0}
\definecolor{lcol}{RGB}{204,102,0}
\definecolor{bcol}{RGB}{0,0,0}
\definecolor{ocol}{RGB}{167,167,167}

\newcommand{\threeaxes}[8]{
	\tikzmath{
    	\Lxr=#3;\Lxl=#3;\Lyt=#4;\Lyb=#4;\Zt= #5;\Zb= #6;\Xang=#7;\Yang=#8;
    	\Xxr= cos(\Xang)*\Lxr; \Xyr=-sin(\Xang)*\Lxr;
    	\Xxl=-cos(\Xang)*\Lxl; \Xyl= sin(\Xang)*\Lxl;
        \Yxt= sin(\Yang)*\Lyt; \Yyt= cos(\Yang)*\Lyt;
        \Yxb=-sin(\Yang)*\Lyb; \Yyb=-cos(\Yang)*\Lyb;
        \zzz=0;
    }
    \begin{scope}[shift={(#1,#2)},rotate=0]
        \ifx\Zt\zzz\else  \draw[black,->] (0,0) -- +(0, \Zt);    \fi
        \ifx\Zb\zzz\else  \draw[black,->] (0,0) -- +(0,-\Zb);    \fi
        \ifx\Lxr\zzz\else \draw[black,->] (0,0) -- +(\Xxr,\Xyr); \fi
        \ifx\Lxl\zzz\else \draw[black,->] (0,0) -- +(\Xxl,\Xyl); \fi
        \ifx\Lyt\zzz\else \draw[black,->] (0,0) -- +(\Yxt,\Yyt); \fi
        \ifx\Lyb\zzz\else \draw[black,->] (0,0) -- +(\Yxb,\Yyb); \fi
	\end{scope}
}
\newcommand{\threeaxeslabelx}[9]{
	\tikzmath{
    	\Lxr=#3;\Lxl=#3;\Lyt=#4;\Lyb=#4;\Zt= #5;\Zb= #6;\Xang=#7;\Yang=#8;
    	\Xxr= cos(\Xang)*\Lxr; \Xyr=-sin(\Xang)*\Lxr;
    	\Xxl=-cos(\Xang)*\Lxl; \Xyl= sin(\Xang)*\Lxl;
        \Yxt= sin(\Yang)*\Lyt; \Yyt= cos(\Yang)*\Lyt;
        \Yxb=-sin(\Yang)*\Lyb; \Yyb=-cos(\Yang)*\Lyb;
        \zzz=0;
    }
    \begin{scope}[shift={(#1,#2)},rotate=0]
        \ifx\Lxr\zzz\else \draw (\Xxr,\Xyr) node[anchor=south west] {#9}; \fi
	\end{scope}
}
\newcommand{\threeaxeslabely}[9]{
	\tikzmath{
    	\Lxr=#3;\Lxl=#3;\Lyt=#4;\Lyb=#4;\Zt= #5;\Zb= #6;\Xang=#7;\Yang=#8;
    	\Xxr= cos(\Xang)*\Lxr; \Xyr=-sin(\Xang)*\Lxr;
    	\Xxl=-cos(\Xang)*\Lxl; \Xyl= sin(\Xang)*\Lxl;
        \Yxt= sin(\Yang)*\Lyt; \Yyt= cos(\Yang)*\Lyt;
        \Yxb=-sin(\Yang)*\Lyb; \Yyb=-cos(\Yang)*\Lyb;
        \zzz=0;
    }
    \begin{scope}[shift={(#1,#2)},rotate=0]
        \ifx\Lyt\zzz\else \draw (\Yxt,\Yyt) node[anchor=south west] {#9}; \fi
	\end{scope}
}
\newcommand{\threeaxeslabelz}[9]{
	\tikzmath{
    	\Lxr=#3;\Lxl=#3;\Lyt=#4;\Lyb=#4;\Zt= #5;\Zb= #6;\Xang=#7;\Yang=#8;
    	\Xxr= cos(\Xang)*\Lxr; \Xyr=-sin(\Xang)*\Lxr;
    	\Xxl=-cos(\Xang)*\Lxl; \Xyl= sin(\Xang)*\Lxl;
        \Yxt= sin(\Yang)*\Lyt; \Yyt= cos(\Yang)*\Lyt;
        \Yxb=-sin(\Yang)*\Lyb; \Yyb=-cos(\Yang)*\Lyb;
        \zzz=0;
    }
    \begin{scope}[shift={(#1,#2)},rotate=0]
        \ifx\Zt\zzz\else  \draw (0,\Zt)     node[anchor=west] {$\;$#9}; \fi
	\end{scope}
}

\newcommand{\isoaxes}[9]{
    \tikzmath{
        \rot=#3;
        \len=#4;
        \ddd=#5;
        \pX=  0; \Xx=cos(\pX)*\len; \Xy=sin(\pX))*\len;
        \pY=120; \Yx=cos(\pY)*\len; \Yy=sin(\pY))*\len;
        \pZ=240; \Zx=cos(\pZ)*\len; \Zy=sin(\pZ))*\len;
        \Xxx=cos(\ddd)*\Xx-sin(\ddd)*\Xy; \Xxy=sin(\ddd)*\Xx+cos(\ddd)*\Xy;
        \Yxx=cos(\ddd)*\Yx-sin(\ddd)*\Yy; \Yxy=sin(\ddd)*\Yx+cos(\ddd)*\Yy;
        \Zxx=cos(\ddd)*\Zx-sin(\ddd)*\Zy; \Zxy=sin(\ddd)*\Zx+cos(\ddd)*\Zy;
    }
    \begin{scope}[shift={(#1,#2)},rotate=\rot]
		\filldraw[black] (0,0) circle (2pt);
        \draw[black,thick,->] (0,0) -- +(\Xx,\Xy);
        \draw[black,thick,->] (0,0) -- +(\Yx,\Yy);
        \draw[black,thick,->] (0,0) -- +(\Zx,\Zy);
        \draw (0.1,0.6)   node[rotate=0,anchor=center] {#6};
        \draw (\Xxx,\Xxy) node[rotate=0,anchor=center] {#7};
        \draw (\Yxx,\Yxy) node[rotate=0,anchor=center] {#8};
        \draw (\Zxx,\Zxy) node[rotate=0,anchor=center] {#9};
    \end{scope}
}

\newcommand{\pane}[7]{
	\tikzmath{
    	\rot=#3;
    	\width=#4;
        \height=#5;
        \corner=#6;
    	\px1= 0.5*\width-\corner; \py1=-0.5*\height;
        \px2= 0.5*\width;         \py2=-0.5*\height+\corner;
        \px3= 0.5*\width;         \py3= 0.5*\height-\corner;
        \px4= 0.5*\width-\corner; \py4= 0.5*\height;
        \px5=-0.5*\width+\corner; \py5= 0.5*\height;
        \px6=-0.5*\width;         \py6= 0.5*\height-\corner;
        \px7=-0.5*\width;         \py7=-0.5*\height+\corner;
        \px8=-0.5*\width+\corner; \py8=-0.5*\height;
    }
	\begin{scope}[shift={(#1,#2)},rotate=\rot]
    	\filldraw[{#7}]
        (\px1,\py1) to[out=    0,in=  -90] (\px2,\py2) --
        (\px3,\py3) to[out=   90,in=    0] (\px4,\py4) --
        (\px5,\py5) to[out= -180,in=   90] (\px6,\py6) --
        (\px7,\py7) to[out=  -90,in= -180] (\px8,\py8) -- cycle;
    \end{scope}
}
\newcommand{\cpane}[7]{
	\tikzmath{
		\pLTx=#1;
		\pRBx=#2;
		\pLTy=#3;
		\pRBy=#4;
	}
	\pane{0.5*\pLTx+0.5*\pRBx}{0.5*\pLTy+0.5*\pRBy}{#5}{\pRBx-\pLTx}{\pLTy-\pRBy}{#6}{#7}
}

\newcommand{\coneback}[7]{
	\tikzmath{\rot=#3;
              \length=#4; \radius=\length*tan(0.5*#5); \depth=#6;
              \sx =  cos(\rot)*#1 + sin(\rot)*#2;
              \sy = -sin(\rot)*#1 + cos(\rot)*#2;
    }
    \begin{scope}[shift={(\sx,\sy)},transform canvas={rotate=\rot}]
	    \draw[{#7}] (\radius,-\length) arc(360:180: {\radius} and {-\depth});
    \end{scope}
}

\newcommand{\cone}[7]{
	\tikzmath{\rot=#3;
              \length=#4; \radius=\length*tan(0.5*#5); \depth=#6;
              \sx =  cos(\rot)*#1 + sin(\rot)*#2;
              \sy = -sin(\rot)*#1 + cos(\rot)*#2;
	}
    \begin{scope}[shift={(\sx,\sy)},transform canvas={rotate=\rot}]
    	\filldraw[{#7}] (0,0) -- (\radius,-\length) arc(360:180: {\radius} and {\depth}) -- cycle;
    \end{scope}
}

\newcommand{\rocket}[6]{
	\tikzmath{\rot=#3;
    		  \length=#4;
              \throttle=(\length/0.8)*0.6*#5;
              \gimbalangle=#6;
    		  \LL = \length; \RR = \LL/8;
              \HH = \LL/5;   \WW = \LL/16;
              \rrrr = \LL/16;  \ww = \LL/16;
              \HG = \LL/8;   \RG = \LL/11; \WG = \LL/16;
              \sx =  cos(\rot)*#1 + sin(\rot)*#2;
              \sy = -sin(\rot)*#1 + cos(\rot)*#2;
    }
    \begin{scope}[shift={(\sx,\sy)},transform canvas={rotate=\rot}]
		\begin{scope}[shift={(0,-0.5*\LL)},rotate=\gimbalangle]
          	\fill[fill=orange!100,shading=axis,shading angle=90,left color=orange!100,right color=orange!25]
            	(-\RG,-\HG) -- (0,-\HG-\throttle) -- (\RG,-\HG) -- cycle;
        	\filldraw[color=black!100,fill=black!10,shading=axis,shading angle=90,left color=black!30,right color=black!0]
        	(0,0) -- (\RG,-\HG) arc(360:180: {\RG} and {\WG}) -- cycle;
		\end{scope}
        
    	\filldraw[color=black!100,fill=black!10,shading=ball,shading angle=90,left color=black!30,right color=black!0]
        	(-\RR,-0.5*\LL) arc(180:360: {\RR} and {\WW}) -- (\RR,0.5*\LL) -- (\RR,0.5*\LL) arc(360:180: {\RR} and {\WW}) -- (-\RR,-0.5*\LL) -- cycle;
        
    	\filldraw[color=black!100,fill=black!10,shading=axis,shading angle=90,left color=black!30,right color=black!0]
        	(\RR,0.5*\LL) arc(360:180: {\RR} and {\WW}) -- (-\rrrr,0.5*\LL+\HH) arc(-180:0: {\rrrr} and {-\ww}) -- cycle;
    \end{scope}
}

\title{Successive Convexification \\ for 6-DoF Powered Descent Guidance \\ with Compound State-Triggered Constraints}

\author{
  Michael Szmuk,
  Taylor P. Reynolds,
  Beh\c cet\ A\c c\i kme\c se,
  Mehran Mesbahi,\\
  {\normalsize\itshape University of Washington, Seattle, WA 98195-2400, USA} \\
  and \\
  John M. Carson III \\
  {\normalsize\itshape NASA Johnson Space Center, Houston, TX 77058, USA}
}

\begin{document}
\maketitle


\begin{abstract}

This paper introduces a continuous formulation for \textit{compound} state-triggered constraints, which are generalizations of the recently introduced state-triggered constraints. State-triggered constraints are different from ordinary constraints found in optimal control in that they use a state-dependent \textit{trigger condition} to enable or disable a \textit{constraint condition}, and can be expressed as continuous functions that are readily handled by successive convexification. Compound state-triggered constraints go a step further, giving designers the ability to compose trigger and constraint conditions using Boolean \textit{and} and \textit{or} operations. Simulations of the 6-degree-of-freedom (DoF) powered descent guidance problem obtained using successive convexification are presented to illustrate the utility of state-triggered and compound state-triggered constraints. The examples employ a velocity-triggered angle of attack constraint to alleviate aerodynamic loads, and a collision avoidance constraint to avoid large geological formations. In particular, the velocity-triggered angle of attack constraint demonstrates the ability of state-triggered constraints to introduce new constraint phases to the solution without resorting to combinatorial techniques.

\end{abstract}


\section*{Nomenclature}
\begin{multicols}{2}
	\begin{tabbing}
	\hspace*{1.6cm}\= \kill
	DoF \> Degrees-of-Freedom \\ 
	DCM \> Direction cosine matrix \\
	STC \> State-Triggered Constraint \\
	cSTC \> Continuous State-Triggered Constraint \\
	$\bvec{e}_i$ \> Unit vector pointing along $i^{th}$-axis \\
	$\inertial$ \> Subscript used to denote the inertial frame \\
	$\body$ \> Subscript used to denote the body frame \\
	$\iframe$ \> The inertial UEN frame \\
	$\bframe$ \> The body frame \\
	$t$ \> Time \\
	$\tin$ \> Initial time epoch \\
	$\tig$ \> Ignition time epoch \\
	$\tf$ \> Final time epoch \\
	$\tc$ \> Coast time \\
	$\tcmax$ \> Max allowable coast time \\
	$\tb$ \> Burn time \\
	$\mdry$ \> Dry mass of the vehicle \\
	$\mig$ \> Mass of vehicle at~$\tig$ \\
	$\gI$ \> Gravity vector \\
	$\m$ \> Vehicle mass \\
	$\rI$ \> Inertial position of the vehicle \\
	$\vI$ \> Inertial velocity of the vehicle \\
	$\qIB$ \> Unit quaternion relating~$\iframe$ to~$\bframe$ \\
	$\omegaB$ \> Angular velocity vector \\
	$\TB$ \> Commanded thrust vector \\
	$\cIB$ \> DCM rotating from~$\iframe$ to~$\bframe$ \\
	$\cBI$ \> DCM rotating from~$\bframe$ to~$\iframe$ \\
	$\glideslope$ \> Glide-slope cone constraint angle \\
	$\tiltmax$ \> Maximum allowable tilt angle \\
	$\omegamax$ \> Maximum allowable angular rate \\
	$\gimbalmax$ \> Maximum allowable gimbal angle \\
	$\Tmin$ \> Minimum allowable thrust magnitude \\
	$\Tmax$ \> Maximum allowable thrust magnitude \\
	\end{tabbing}
\end{multicols}

\newpage


\section{Introduction} \label{sec:introduction}

Real-time optimal powered descent guidance algorithms that enable precision landing are key to maximizing the probability of landing success. Such algorithms allow the vehicle to cope with a more adverse set of model uncertainties, environmental disturbances, and unexpected obstacles as the vehicle approaches the landing site. This increase in robustness can be used to attempt landings in more challenging and scientifically interesting areas, or instead can be reallocated by exchanging excessive propellant mass for an increased payload mass fraction. In recent years, optimal powered descent guidance technology has played a key role in the robust recovery of commercial vertical-takeoff-vertical-landing reusable launch vehicles~\cite{larsNAE}, and is at the forefront of the drive to reduce launch costs for the foreseeable future.

Solving the powered descent guidance problem quickly and reliably is necessary because the vehicle has a limited amount of propellant, and must react quickly to deviations and conditions observed close to the ground. Doing so is challenging due to the nonlinear nature of the dynamics, the non-convex nature of the state and control constraints, and the free-final-time nature of the problem.

Work on powered descent guidance began during the Apollo program, with~\cite{Meditch1964,Lawden1963,Marec1979} approaching the problem using optimal control theory and calculus of variations. However, these methods were not incorporated into the Apollo flight code since the polynomial-based guidance methods in use at the time were deemed sufficiently optimal, and were far simpler to design and implement~\cite{klumpp}. After the Apollo program, research continued in search of analytical solutions to the 3-degree-of-freedom (DoF) landing problem~\cite{Kornhauser1972,Breakwell1975,Azimov1996,dsouza97}.

In the early 2000s, interest in the powered descent landing problem renewed, this time with a focus on robotic missions to Mars. A number of works using direct methods were published, with~\cite{Topcu2005,topcu_pdg} using numerical simulations to demonstrate theoretical results on the 3-DoF problem, and with~\cite{Acikmese2005,ploenaiaa06,behcetjgcd07} culminating in the use of Pontryagin's maximum principle to losslessly convexify the 3-DoF problem. These works were later generalized in~\cite{behcetaut11,harris_acc,lars_12scl,behcet_aut11,larssys12,pointing2013,matt_aut1}, and were demonstrated in a sequence of flight experiments in the early 2010s~\cite{dueri2016customized,Scharf14,gfold,gfold_ieee_aero_14}.

In 2015, a dual-quaternion-based approach was proposed in which a 6-DoF line-of-sight constraint was convexified \cite{mesbahi_6dof,lee2016constrained}. This method was inherently equipped to handle 6-DoF motion, but relied on piecewise-affine approximations to deal with the nonlinear dynamics. As such, the solution degraded in accuracy with coarser temporal discretizations.

While the convexification techniques listed above represent significant advancements in powered descent guidance, they remain applicable to a relatively limited class of optimal control problems. On the other hand, sequential convex programming (SCP) and successive convexification techniques offer a framework for solving more general non-convex optimal control problems, at the expense of more computational complexity and weakened optimality and convergence guarantees~\cite{jordithesis,liu2014solving,liu2015entry,SCvx_cdc16,SCvx_2017arXiv,SCvx_2018arXiv,wang2016constrained,liu_new}. These methods work by converting the non-convex problem into a sequence of convex subproblems that locally approximate the non-convexities to first-order, and do not carry the approximation to second order as in Sequential Quadratic Programming (SQP) methodologies~\cite{Boggs1995}.

The work presented herein uses the successive convexification and discretization framework developed in our earlier work~\cite{szmuk2016successive,szmuk2017successive,szmuk2018successive}. We also draw heavily from our recent work on a free-final-time 6-DoF powered descent guidance problem with state-triggered constraints~\cite{SzmukReynolds2018}. State-triggered constraints represent a novel formulation that, to the best of our knowledge, has not been introduced in the optimal control literature. It bears a close resemblance to complementarity constraints~\cite{Cottle1992,Heemels2000}, which use continuous variables to model discrete events. Complementarity constraints can be more efficient than mixed-integer approaches (see Section~2.3 in~\cite{biegler2014}), which suffer from poor complexity and are not well suited for solving the powered descent guidance problem in real-time~\cite{Richards2015}. However, unlike complementarity constraints that make mutually exclusive decisions between two variables, state-triggered constraints are capable of one-way \textit{if}-\textit{then} implications. The primary contribution of this paper is the introduction of a continuous formulation for \textit{compound} state-triggered constraints, which generalizes our existing state-triggered constraint formulation.

In subsequent sections, we use the following notation.~$\real_{+}$ and~$\real_{++}$ respectively denote non-negative and strictly positive reals,~$\spd{n}$ denotes an~$n\times n$ symmetric positive-definite matrix,~$\mathcal{S}^3$ denotes the unit 3-sphere, and the superscript~$c$ denotes the complement of a set. The symbols $\dotprod{}{}$ and~$\crossprod{}{}$ respectively represent the vector dot and cross products,~$\bvec{a}\perp\bvec{b}$ denotes the orthogonality constraint~$\dotprod{a}{b}=0$, and~$\bvec{e}_{i}$ denotes the unit vector with~$1$ located at the~$\ith{i}{th}$ element. The subscripts~$\inertial$ and~$\body$ respectively denote quantities expressed in inertial frame~($\iframe$) and body frame~($\bframe$) coordinates. The remainder of this paper is organized as follows. In~\sref{sec:stcs}, we discuss state-triggered constraints and the main contribution of this paper; in~\sref{sec:prob_state}, we give an overview of successive convexification; in~\sref{sec:numerical_results}, we present simulation results; and in~\sref{sec:conclusion}, we conclude the paper.

\newpage


\section{Continuous State-Triggered Constraints} \label{sec:stcs}

This section is organized as follows. In~\sref{sec2:stc}, we review the state-triggered constraint (STC) formulation introduced in Section II.F of~\cite{SzmukReynolds2018}. Specifically,~\sref{sec2:stc_form_def} formally defines STCs,~\sref{sec2:stc_cont_form} outlines the continuous formulation of STCs,~\sref{sec2:stc_ex} provides an example application for STCs, and~\sref{sec2:stc_alt_forms} discusses alternative formulations. In~\sref{sec2:cmpstc}, we follow the same outline as in~\sref{sec2:stc} to introduce an extension of state-triggered constraints called \textit{compound state-triggered constraints}.


\subsection{State-Triggered Constraints} \label{sec2:stc}

To better understand state-triggered constraints, let us first consider the most common type of constraints encountered in the optimal control literature: \textit{temporally-scheduled constraints}. Temporally-scheduled constraints are enforced over predefined time intervals that are determined during problem formulation, and that are not altered while the problem is being solved. For example, such constraints are often be enforced over the entire time horizon of an optimal control problem.

Unlike a temporally-schedule constraint, a state-triggered constraint is enforced only when a state-dependent condition is satisfied, and functions as an \textit{if}-statement conditioned on the solution variable it constrains. An optimal control problem containing a state-triggered constraint determines its solution variable with a simultaneous understanding of how the constraint affects the optimization, and of how the optimization enables or disables the constraint.

\subsubsection{Formal Definition}  \label{sec2:stc_form_def}

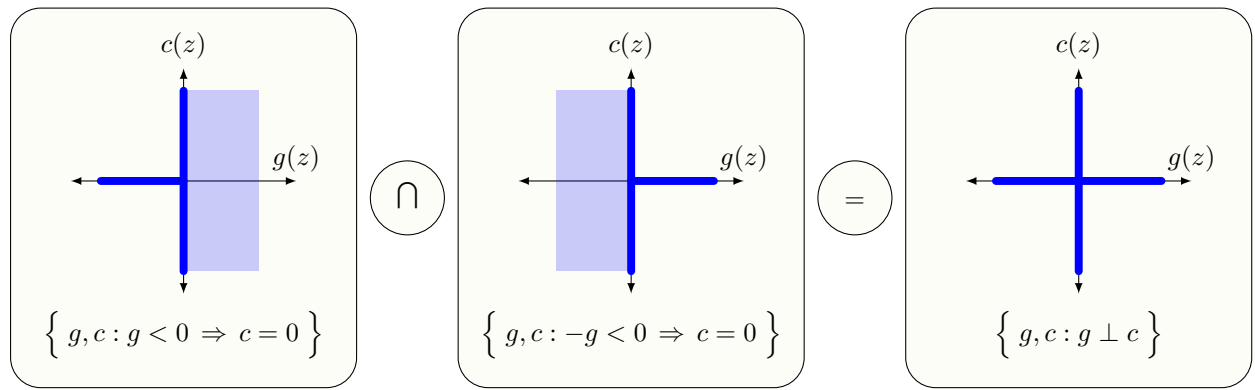
\begin{figure}[b!]
	\centering
  	\begin{tikzpicture}
    \tikzmath{\lx=1.5; \ly=1.5; \dy=0.1;
              \cx1= 5.95;   \cy1=0; 
              \cx2= 0.0;    \cy2=0; 
              \cx3= 2*\cx1; \cy3=0; 
              \llx=1.0; \lly=1.0;
              \lllx=\llx+0.1; \llly=\lly+0.2;
              \hx=0.5*(16.5-\cx3); \hy1=2.3; \hy2=2.75; \hc=0.4;
              \ccy=\cy1-\hy2+0.5*(\hy1+\hy2);}
    
    \filldraw[color=black,fill=beige!20]
    (\cx2-\hx,\cy2-\hy2+\hc) to[out=-90,in=-180] (\cx2-\hx+\hc,\cy2-\hy2) -- 
    (\cx2+\hx-\hc,\cy2-\hy2) to[out=0,in=-90] (\cx2+\hx,\cy2-\hy2+\hc) --
    (\cx2+\hx,\cy2+\hy1-\hc) to[out=90,in=0] (\cx2+\hx-\hc,\cy2+\hy1) --
    (\cx2-\hx+\hc,\cy2+\hy1) to[out=-180,in=90] (\cx2-\hx,\cy2+\hy1-\hc) -- cycle;

    \filldraw[color=black,fill=beige!20]
    (\cx1-\hx,\cy1-\hy2+\hc) to[out=-90,in=-180] (\cx1-\hx+\hc,\cy1-\hy2) -- 
    (\cx1+\hx-\hc,\cy1-\hy2) to[out=0,in=-90] (\cx1+\hx,\cy1-\hy2+\hc) --
    (\cx1+\hx,\cy1+\hy1-\hc) to[out=90,in=0] (\cx1+\hx-\hc,\cy1+\hy1) --
    (\cx1-\hx+\hc,\cy1+\hy1) to[out=-180,in=90] (\cx1-\hx,\cy1+\hy1-\hc) -- cycle;

    \filldraw[color=black,fill=beige!20]
    (\cx3-\hx,\cy3-\hy2+\hc) to[out=-90,in=-180] (\cx3-\hx+\hc,\cy3-\hy2) -- 
    (\cx3+\hx-\hc,\cy3-\hy2) to[out=0,in=-90] (\cx3+\hx,\cy3-\hy2+\hc) --
    (\cx3+\hx,\cy3+\hy1-\hc) to[out=90,in=0] (\cx3+\hx-\hc,\cy3+\hy1) --
    (\cx3-\hx+\hc,\cy3+\hy1) to[out=-180,in=90] (\cx3-\hx,\cy3+\hy1-\hc) -- cycle;
	
    \filldraw[color=black,fill=beige!20] (0.5*\cx1+0.5*\cx2,\ccy) circle (14pt);
    \draw (0.5*\cx1+0.5*\cx2,\ccy-0.26) node[anchor=south] {$\bigcap$};
    
    \filldraw[color=black,fill=beige!20] (0.5*\cx1+0.5*\cx3,\ccy) circle (14pt);
    \draw (0.5*\cx1+0.5*\cx3,\ccy-0.24) node[anchor=south] {$=$};
    
    \draw[black,->] (\cx2,\cy2) -- +(\lx,0);
    \draw[black,->] (\cx2,\cy2) -- +(-\lx,0);
    \draw[black,->] (\cx2,\cy2) -- +(0,\ly);
    \draw[black,->] (\cx2,\cy2) -- +(0,-\ly);
    \draw (\cx2+\lx,0) node[anchor=south] {$\stcg(\stcx)$};
    \draw (\cx2,\cy2+\ly) node[anchor=south] {$\stcf(\stcx)$};
    \draw (\cx2,\cy2-\ly-\dy) node[anchor=north] {$\Big\{\;\stcg,\stcf: \stcg < 0 \,\Rightarrow\, \stcf=0 \;\Big\}$};
    \draw[blue,line width=3pt,line cap=round] (\cx2,\cy2-\llly) -- (\cx2,\cy2+\llly);
    \draw[blue,line width=3pt,line cap=round] (\cx2,\cy2) -- (\cx2-\lllx,\cy2);
    \fill[blue,opacity=0.2] (\cx2,\cy2+\llly)
                         -- (\cx2+\llx,\cy2+\llly)
                         -- (\cx2+\llx,\cy2-\llly)
                         -- (\cx2,\cy2-\llly) -- cycle;
    
    \draw[black,->] (\cx1,\cy1) -- +(\lx,0);
    \draw[black,->] (\cx1,\cy1) -- +(-\lx,0);
    \draw[black,->] (\cx1,\cy1) -- +(0,\ly);
    \draw[black,->] (\cx1,\cy1) -- +(0,-\ly);
    \draw (\cx1+\lx,0) node[anchor=south] {$\stcg(\stcx)$};
    \draw (\cx1,\cy1+\ly) node[anchor=south] {$\stcf(\stcx)$};
    \draw (\cx1,\cy1-\ly-\dy) node[anchor=north] {$\Big\{\;\stcg,\stcf: -\stcg < 0 \,\Rightarrow\, \stcf=0\;\Big\}$};
    \draw[blue,line width=3pt,line cap=round] (\cx1,\cy1-\llly) -- (\cx1,\cy1+\llly);
    \draw[blue,line width=3pt,line cap=round] (\cx1,\cy1) -- (\cx1+\lllx,\cy1);
    \fill[blue,opacity=0.2] (\cx1,\cy1+\llly)
                         -- (\cx1-\llx,\cy1+\llly)
                         -- (\cx1-\llx,\cy1-\llly)
                         -- (\cx1,\cy1-\llly) -- cycle;

    \draw[black,->] (\cx3,\cy3) -- +(\lx,0);
    \draw[black,->] (\cx3,\cy3) -- +(-\lx,0);
    \draw[black,->] (\cx3,\cy3) -- +(0,\ly);
    \draw[black,->] (\cx3,\cy3) -- +(0,-\ly);
    \draw (\cx3+\lx,0) node[anchor=south] {$\stcg(\stcx)$};
    \draw (\cx3,\cy3+\ly) node[anchor=south] {$\stcf(\stcx)$};
    \draw (\cx3,\cy3-\ly-\dy) node[anchor=north] {$\Big\{\;\stcg,\stcf: \stcg\perp\stcf \;\Big\}$};
    \draw[blue,line width=3pt,line cap=round] (\cx3,\cy3-\llly) -- (\cx3,\cy3+\llly);
    \draw[blue,line width=3pt,line cap=round] (\cx3-\lllx,\cy3) -- (\cx3+\lllx,\cy3);
\end{tikzpicture}
  	\vspace{-0.2cm}
    \caption{The feasible sets of STCs (left and center panes) and orthogonality constraints (right pane). The feasible sets of the STCs can be intersected in order to recover the feasible set of an orthogonality constraint, and thus capture a broader set of logic in the~$\stcg$-$\stcf$ plane.}
  	\label{fig:stc_vs_oc}
\end{figure}

A state-triggered constraint is comprised of a \textit{trigger condition} and a \textit{constraint condition} that are arranged such that satisfaction of the former implies satisfaction of the latter. More formally, the logical implication of a state-triggered constraint is given by
\begin{equation} \label{eq:stc_formal_def}
    \stcg(\stcx)<0 \;\Rightarrow\; \stcf(\stcx)=0,
\end{equation}
where~$\stcx\in\real^\stcnx$ represents the solution variable of the parent optimization problem;~$\stcg(\cdot) : \real^\stcnx \rightarrow \real$ and~$\stcf(\cdot) : \real^\stcnx \rightarrow \real$ are piecewise continuously differentiable functions that are termed the \textit{trigger function} and \textit{constraint function}, respectively; and~$\stcg(\stcx)<0$ and~$\stcf(\stcx) = 0$ are the aforementioned trigger and constraint conditions.

The feasible set corresponding to~\eqref{eq:stc_formal_def} is shown in the leftmost pane of Figure~\ref{fig:stc_vs_oc}. The figure illustrates that simultaneously enforcing the state-triggered constraints~$\stcg(\stcx)<0\,\Rightarrow\,\stcf(\stcx)=0$ and~$-\stcg(\stcx)<0 \,\Rightarrow\,\stcf(\stcx)=0$ yields the feasible set corresponding to the orthogonality constraint~$\stcg(\stcx)\perp\stcf(\stcx)$. Further, intersecting these feasible sets with the set~$\{\,\stcg,\stcf:\stcg\geq 0\,,\,\stcf\geq 0\,\}$ recovers the feasible set of the well-known complementarity constraint~$0 \leq\stcg(\stcx)\perp\stcf(\stcx)\geq 0$~\cite{Cottle1992}. When viewed in the~$\stcg$-$\stcf$ plane, both orthogonality and complementarity constraints effectively represent bi-directional \textit{if-and-only-if}-statements, whereas state-triggered constraints emulate \textit{if}-statements. Thus, we argue that the state-triggered constraint formulation represents a more general logical building block for the design of guidance algorithms.

\begin{figure}[t!]
    \centering
    \begin{tikzpicture}
    \newcommand{\lw}{0.2mm}
    \newcommand{\llw}{0.5mm}
    \newcommand{\rrrad}{2.00}
    \definecolor{mypurple}{rgb}{0.0,0.6,0.0}
    \definecolor{mygreen}{rgb}{0.0,0.6,0.0}

    \tikzmath{
        \ldd=0.5;
        \lx=16.5;   \ly=7.0;
        \lddtx=0.1; \lddty=0.1;
        \cccx=1.5;  \cccy=0.0;
        \pppx=\cccx-1.50; \pppy=\cccy+0.75;
        \qqqx=\cccx+1.50; \qqqy=\cccy+0.75;
        \rrrx=\cccx+0.00; \rrry=\cccy-0.75;
        \dddr=0.15;
        \aaa0= 25;
        \aaa1=  46.0; \bbb1= 314.5;
        \aaa2=-143;   \bbb2=-217.5;
        \aaa3=-190;   \bbb3=-135;
        \aaa4=-159.5; \bbb4=100;
        \aaa5=-38;    \bbb5=-96;
        \aaa6= 35;
        \aaaq=-20;
        \ppppx=\pppx-(\rrrad-\dddr)*cos(\aaa0);
        \ppppy=\pppy+(\rrrad-\dddr)*sin(\aaa0);
        \qqqqx=\pppx-(\rrrad)*cos(\aaaq);
        \qqqqy=\pppy+(\rrrad)*sin(\aaaq);
        \pppppx=\rrrx+(\rrrad-\dddr)*cos(\aaa6);
        \pppppy=\rrry-(\rrrad-\dddr)*sin(\aaa6);
    }
    
    \pane{0}{0}{0}{\lx}{\ly}{0.4}{color=black,fill=beige!20};
    
    \filldraw [fill=gray!1,draw=black,line width=0.25mm] (-0.5*\lx+\ldd,0.5*\ly-\ldd) rectangle (0.5*\lx-\ldd,-0.5*\ly+\ldd);
    
    \fill [fill=red    ,fill opacity=0.10] (\pppx,\pppy) circle(\rrrad);
    \fill [fill=blue   ,fill opacity=0.10] (\qqqx,\qqqy) circle(\rrrad);
    \fill [fill=mygreen,fill opacity=0.15] (\rrrx,\rrry) circle(\rrrad);
    
    \filldraw [draw=red,line width=\llw,cap=round,fill=red!25]
        ({\pppx+(\rrrad-\dddr)*cos(\aaa1)},
         {\pppy+(\rrrad-\dddr)*sin(\aaa1)}) arc (\aaa1:\bbb1:\rrrad-\dddr) --
        ({\qqqx+(\rrrad+\dddr)*cos(\aaa2)},
         {\qqqy+(\rrrad+\dddr)*sin(\aaa2)}) arc (\aaa2:\bbb2:\rrrad+\dddr) -- cycle;
    
    \filldraw [draw=mypurple,line width=\llw,cap=round,fill=mypurple!25]
        ({\qqqx+(\rrrad-\dddr)*cos(\aaa3)},
         {\qqqy+(\rrrad-\dddr)*sin(\aaa3)}) arc (\aaa3:\bbb3:\rrrad-\dddr) --
        ({\pppx+(\rrrad+\dddr)*cos(\aaa5)},
         {\pppy+(\rrrad+\dddr)*sin(\aaa5)}) arc (\aaa5:\bbb5:\rrrad+\dddr) --
        ({\rrrx+(\rrrad-\dddr)*cos(\aaa4)},
         {\rrry+(\rrrad-\dddr)*sin(\aaa4)}) arc (\aaa4:\bbb4:\rrrad-\dddr) -- cycle;
    
    \draw [color=red    ,line width=\lw] (\pppx,\pppy) circle(\rrrad);
    \draw [color=blue   ,line width=\lw] (\qqqx,\qqqy) circle(\rrrad);
    \draw [color=mygreen,line width=\lw] (\rrrx,\rrry) circle(\rrrad);
    
    \draw [black] (0.5*\lx-\ldd-\lddtx,0.5*\ly-\ldd-\lddty) node[anchor=north east] {$\real^{\nz}$};
    \draw [darkgreen] (\cccx-2.1,\cccy-1.7) node[anchor=center] {$\mathcal{Z}$};
    \draw [darkred] (\cccx-3.8,\cccy+0.8) node[anchor=center] {$\mathcal{G}$};
    \draw [blue] (\cccx+3.5,\cccy-0.3) node[anchor=center] {$\mathcal{C}$};
    
    \filldraw[black] (\ppppx,\ppppy) circle(1pt);
    \draw[black,line width=0.1mm] (-0.5*\lx+\ldd+\lddtx+4.3,0.5*\ly-\ldd-\lddty-1.00) to[out=0,in=180-\aaa0] (\ppppx,\ppppy);
    \draw[black] (-0.5*\lx+\ldd+\lddtx-0.1,0.5*\ly-\ldd-\lddty) node[anchor=north west] {\footnotesize\it \begin{tabular}{l}
        The state-triggered constraint \\ $\stcg(\stcx)<0\,\Rightarrow\,\stcf(\stcx)=0$ ensures \\
        that feasible solutions do \\
        not lie inside~$\mathcal{G}\cap\mathcal{C}^{c}$ \end{tabular}};
    
    
    
    
    \filldraw[black] (\cccx,\cccy+0.5) circle(1pt);
    \draw[black,line width=0.1mm] (0.5*\lx-\ldd-\lddtx-2.05,-0.5*\ly+\ldd+\lddty+3.90) to[out=220,in=-10] (\cccx,\cccy+0.5);
    \draw[black] (0.5*\lx-\ldd-\lddtx,-0.5*\ly+\ldd+\lddty+3.3) node[anchor=south east] {\footnotesize\it \begin{tabular}{r}Set where the \\ state-triggered \\ constraint \\ is active\end{tabular}};
    
    \filldraw[black] (\pppppx,\pppppy) circle(1pt);
    \draw[black,line width=0.1mm] (0.5*\lx-\ldd-\lddtx-3.5,-0.5*\ly+\ldd+\lddty+0.70) to[out=180,in=-\aaa6] (\pppppx,\pppppy);
    \draw[black] (0.5*\lx-\ldd-\lddtx,-0.5*\ly+\ldd+\lddty) node[anchor=south east] {\footnotesize\it \begin{tabular}{r}
        The feasible set \\
        with the state-triggered  \\
        constraint included is~$\mathcal{Z}\cap\big(\mathcal{G}^{c}\cup\mathcal{C}\big)$ \end{tabular}};
    
    \draw [black] (-0.5*\lx+\ldd+\lddtx,-0.5*\ly+\ldd+\lddty) node[anchor=south west] {\footnotesize
        $ \begin{aligned}
            \mathcal{G}&\definedas\big\{\,\stcx\in\real^{\nz}: \stcg(\stcx) < 0\,\big\} \\
            \mathcal{C}&\definedas\big\{\,\stcx\in\real^{\nz}: \stcf(\stcx) = 0\,\big\} \\
            \mathcal{Z}&\definedas\big\{\,\stcx\in\real^{\nz}: \stcx\;\textit{is feasible w.r.t. all other constraints}\,\big\}
        \end{aligned} $};
\end{tikzpicture}
\vspace{0.1cm}
  	\vspace{-0.2cm}
    \caption{A Venn diagram of the sets of solution variables~$\stcx\in\real^{\nz}$ that satisfy the trigger condition ($\mathcal{G}$), constraint condition ($\mathcal{C}$), and all other state and control constraints excluding STCs ($\mathcal{Z}$). State-triggered constraints ensure that feasible solutions satisfy~$z\notin\big(\mathcal{G}\cap\mathcal{C}^{c}\big)$.}
  	\label{fig:stc_venn}
  	\vspace{-0.5cm}
\end{figure}
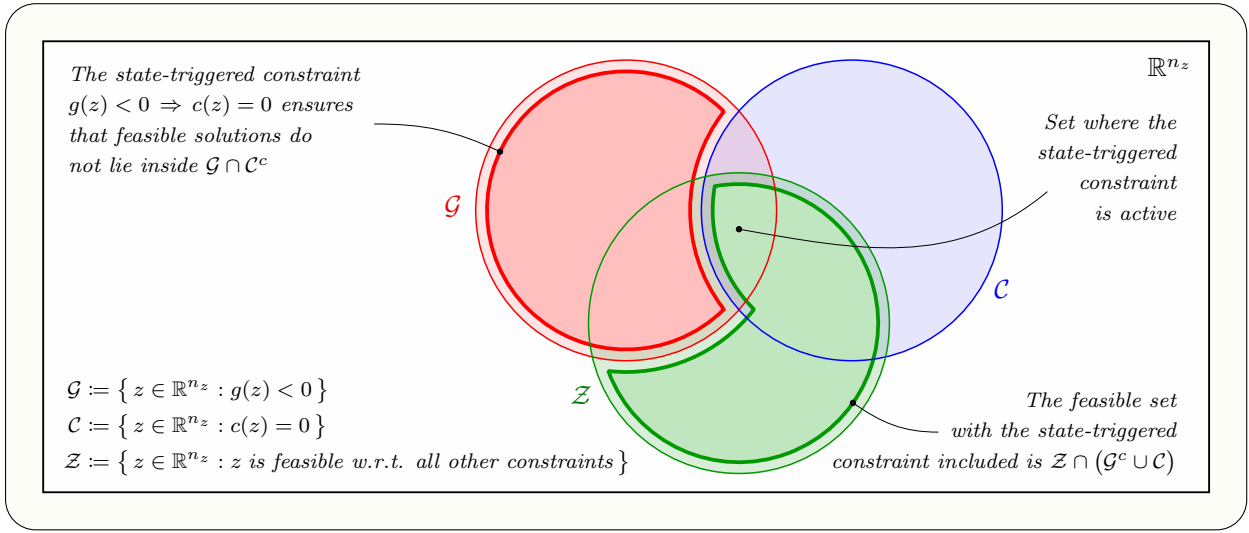

Defining the sets~$\mathcal{G}$,~$\mathcal{C}$, and~$\mathcal{Z}$ as in Figure~\ref{fig:stc_venn}, the statement ``$\stcx\in\mathcal{G}$ implies~$\stcx\in\mathcal{C}\,$'' is typically taken to mean that~$\mathcal{G}\subseteq\mathcal{C}$. However, as seen in Figure~\ref{fig:stc_venn},~$\mathcal{G}$ is not necessarily a subset of~$\mathcal{C}$. In fact, if the trigger and constraint functions are defined such that~$\mathcal{G}\subseteq\mathcal{C}$, then the corresponding state-triggered constraint is trivially satisfied and serves no practical purpose in the context of the optimization problem. However, if~$\stcg(\stcx)$ and~$\stcf(\stcx)$ are defined such that~$\mathcal{G}\not\subseteq\mathcal{C}$, then the state-triggered constraint \textit{removes} the set~$\mathcal{G}\cap\mathcal{C}^{c}$ from the feasible set, thus ensuring that $\mathcal{G}\subseteq\mathcal{C}$ over the remaining feasible set.

Before proceeding, we note that~\eqref{eq:stc_formal_def} can equivalently represent a state-triggered \textit{inequality} constraint, provided that~$\stcx$ is augmented with a non-negative slack variable and that~$\stcg(\stcx)$ and~$\stcf(\stcx)$ are modified accordingly~\cite{BoydConvex}. Thus, we proceed with the understanding that the equality in~\eqref{eq:stc_formal_def} can just as well be replaced with an inequality. Additionally, from De Morgan's Law, satisfaction of~\eqref{eq:stc_formal_def} implies satisfaction of the contrapositive~$\stcf(\stcx)\neq 0\,\Rightarrow\,\stcg(\stcx)\neq 0$. This behavior is exhibited in the simulation examples presented in~\sref{sec:numerical_results}.

\subsubsection{Continuous Formulation} \label{sec2:stc_cont_form}

The continuous formulation corresponding to the state-triggered constraint in~\eqref{eq:stc_formal_def} is given by the following system of equations
\begin{subequations} \label{eq:stc_cont_std}
	\begin{align}
    	\stcg(\stcx) + \stcs &\geq 0, \label{eq:stc_cont_std_a} \\
        \stcs &\geq 0, \label{eq:stc_cont_std_b} \\
        \stcs\cdot\stcf(\stcx) &= 0, \label{eq:stc_cont_std_c}
    \end{align}
\end{subequations}
where~$\stcs\in\real_{+}$ is an auxiliary variable that enables or disables the constraint condition. We refer to~\eqref{eq:stc_cont_std} as the \textit{original formulation}, and note that it possesses an ambiguity in~$\stcs$ given a~$\stcg(\stcx)$ and~$\stcf(\stcx)$. This form can be written more compactly and without the ambiguity in~$\stcs$ using the \textit{projected formulation} given by
\begin{equation} \label{eq:stc_cont_proj}
    \stch(\stcx)\definedas -\min\big(\stcg(\stcx),0\big)\cdot\stcf(\stcx) = 0.
\end{equation}
Note that the feasible sets of~\eqref{eq:stc_cont_std} and~\eqref{eq:stc_cont_proj} are identical to that of~\eqref{eq:stc_formal_def} shown in the leftmost pane of Figure~\ref{fig:stc_vs_oc}. The \textit{continuous state-triggered constraints} (cSTCs) presented above are therefore logically equivalent to the state-triggered constraints detailed in~\sref{sec2:stc_form_def}, and can be readily incorporated into a continuous optimization framework. For more details, we refer the reader to Section~II.F in~\cite{SzmukReynolds2018}.

\subsubsection{Example Application} \label{sec2:stc_ex}

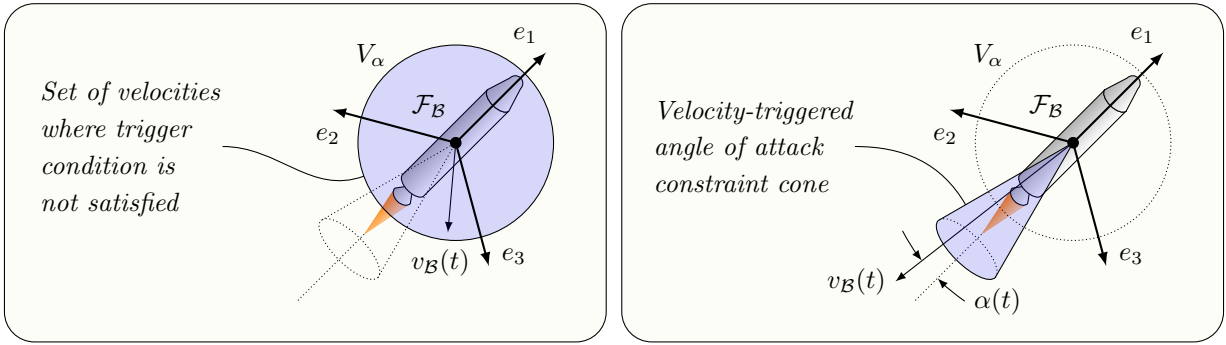
\begin{figure}[t!]
    \begin{subfigure}{0.49\textwidth}
    	\centering
      	\begin{tikzpicture}
    \tikzmath{\aa=-45;
    		  \cx=2; \cy=0.8;
              \Lcone=2.00; \Acone=30; \Dcone=0.2;
              \labx=-0.8; \laby=0.75;
              \LL=3.0; \aaa=-7; \fact=0.85; \da=10;
              \lx=sin(\aa)*\LL; \ly=-cos(\aa)*\LL;
              \llx=sin(-5)*1.2; \lly=-cos(-5)*1.2;
              \RRRR=1.3; \aaaa=200;
              \lllx=\cx+\RRRR*cos(\aaaa); \llly=\cy+\RRRR*sin(\aaaa);
    }
    
    \pane{0}{0.40}{0}{8.0}{4.5}{0.4}{color=black,fill=beige!20}
    
	\draw[black,densely dotted] (\cx,\cy) -- +(\lx,\ly);
    \coneback{\cx}{\cy}{\aa}{\Lcone}{\Acone}{\Dcone}{color=black!100,densely dotted}
    
    \rocket{\cx}{\cy}{\aa}{1.6}{0.6}{0}
    
    \draw[black!100,->] (\cx,\cy) -- +(\llx,\lly);
    \draw (\cx+\llx-0.1,\cy+\lly-0.1) node[anchor=north] {$\vB(t)$};
    
    \cone{\cx}{\cy}{\aa}{\Lcone}{\Acone}{\Dcone}{draw=black!100,fill=blue!100,fill opacity=0.00,densely dotted};
    
    \draw (\cx-0.8,\cy+0.9) node[anchor=south east] {$\Vaoa$};
    \filldraw[draw=black!100,fill=blue!100,fill opacity=0.15] (\cx,\cy) ellipse ({\RRRR} and {\RRRR});
    \draw (\lllx,\llly) to[out=\aaaa,in=0] (\labx-0.3,\laby);
    \draw (\labx,\laby) node[anchor=east] {\it\begin{tabular}{l}
        Set of velocities \\
        where trigger \\
        condition is \\
        not satisfied \end{tabular}};
    
    \isoaxes{\cx}{\cy}{90+\aa}{1.7}{12}{$\bframe$}{$\bvec{e}_1$}{$\bvec{e}_2$}{$\bvec{e}_3$}
    
\end{tikzpicture}
      	\vspace{-0.2cm}
      	\caption{Disabled Constraint Condition}
      	\label{fig:stc_aoa_disabled}
    \end{subfigure}
    \begin{subfigure}{0.49\textwidth}
        \centering
        \begin{tikzpicture}
    \tikzmath{\aa=-45;
    		  \cx=2; \cy=0.8;
              \Lcone=2.00; \Acone=30; \Dcone=0.2;
              \labx=-0.6; \laby=0.75;
              \LL=3.0; \aaa=-7; \fact=0.85; \da=10;
              \lx=sin(\aa)*\LL; \ly=-cos(\aa)*\LL;
              \llx=sin(\aa+\aaa)*\LL; \lly=-cos(\aa+\aaa)*\LL;
              \RRRR=1.3; \aaaa=200;
              \lllx=\cx+sin(\aa-0.5*\Acone)*0.55*\LL; \llly=\cy-cos(\aa-0.5*\Acone)*0.55*\LL;
    }
    
    \pane{0}{0.40}{0}{8.0}{4.5}{0.4}{color=black,fill=beige!20}
    
	\draw[black,densely dotted] (\cx,\cy) -- +(\lx,\ly);
    \coneback{\cx}{\cy}{\aa}{\Lcone}{\Acone}{\Dcone}{color=black!100,densely dotted}
    
    \rocket{\cx}{\cy}{\aa}{1.6}{0.6}{0}
    
    \draw[black!100,->] (\cx,\cy) -- +(\llx,\lly);
    \draw (\cx+\llx,\cy+\lly) node[anchor=east] {$\vB(t)$};
    
    \cone{\cx}{\cy}{\aa}{\Lcone}{\Acone}{\Dcone}{draw=black!100,fill=blue!100,fill opacity=0.15};
    
    \draw (\cx-0.8,\cy+0.9) node[anchor=south east] {$\Vaoa$};
    \filldraw[draw=black!100,fill=blue!100,fill opacity=0.00,densely dotted] (\cx,\cy) ellipse ({\RRRR} and {\RRRR});
    \draw (\lllx,\llly) to[out=180-0.5*\Acone+\aa,in=0] (\labx-0.3,\laby);
    \draw (\labx,\laby) node[anchor=east] {\it\begin{tabular}{l}
        Velocity-triggered \\
        angle of attack \\
        constraint cone \end{tabular}};
    
    \isoaxes{\cx}{\cy}{90+\aa}{1.7}{12}{$\bframe$}{$\bvec{e}_1$}{$\bvec{e}_2$}{$\bvec{e}_3$}
    
    \draw[color=black!100,<-] (\cx+\fact*\lx,\cy+\fact*\ly) arc(-90+\aa:-90+\aa+\da:2*\fact*\LL) node[anchor=west] {$\aoa(t)$};
    \draw[color=black!100,<-] (\cx+\fact*\llx,\cy+\fact*\lly) arc(-90+\aa+\aaa:-90+\aa+\aaa-\da:2*\fact*\LL);
\end{tikzpicture}
      	\vspace{-0.2cm}
        \caption{Enabled Constraint Condition}
        \label{fig:stc_aoa_enabled}
    \end{subfigure}
  	\caption{An example powered descent guidance application where a velocity-triggered angle of attack constraint is used to alleviate aerodynamic loads experienced by the vehicle. The set of velocities that do not satisfy the trigger condition (and for which the constraint condition is not enforced) is represented by the blue sphere in (a). When the trigger condition is satisfied, the velocities are constrained to the blue cone shown in (b).}
  	\label{fig:stc_ex}
  	\vspace{-0.5cm}
\end{figure}

The continuous state-triggered constraint formulation outlined in~\sref{sec2:stc_cont_form} allows a richer set of constraints to be incorporated into a continuous optimization framework. For example, consider the problem of aerodynamic load alleviation during an atmospheric powered descent guidance maneuver. Imposing a single angle of attack constraint over the entire length of the trajectory may be overly conservative, and may eliminate large swaths of otherwise feasible solutions. Hence, we state the objective of an aerodynamic load alleviation constraint as follows: limit the angle of attack of the vehicle when the dynamic pressure is high, and leave the angle of attack unconstrained when the dynamic pressure is negligible.

Assuming constant atmospheric density, this constraint may be interpreted as a velocity-triggered angle of attack constraint, and is stated formally as the following state-triggered constraint
\begin{equation} \label{eq:stc_ex_aoa_orig}
	\twonorm{\vB(t)} > \Vaoa \;\Rightarrow\; -\dotprod{\bvec{e}_1}{\vB(t)} \geq \cos{\aoamax}\twonorm{\vB(t)},
\end{equation}
where~$\Vaoa\in\real_{++}$ is a speed above which the angle of attack is limited to~$\aoa(t) \in\intee{0}{\aoamax}$, and~$\vB(t)\in\real^3$ is the velocity expressed in~$\bframe$ coordinates. Figure~\ref{fig:stc_aoa_disabled} shows the set of feasible velocities when the trigger condition \textit{is not} satisfied and the constraint condition is \textit{disabled}. Figure~\ref{fig:stc_aoa_enabled} shows the set of feasible velocities when the trigger condition \textit{is} satisfied and the constraint condition is \textit{enabled}. Using the projected formulation in~\eqref{eq:stc_cont_proj}, the STC in~\eqref{eq:stc_ex_aoa_orig} is converted into the following cSTC:
\begin{subequations} \label{eq:stc_ex_aoa_imp}
	\begin{align}
		\haoa{\vI(t)}{\qIB(t)} &\definedas -\min\Big(\gaoa{\vI(t)},0 \Big) \cdot \faoa{\vI(t)}{\qIB(t)} \leq 0, \label{eq:stc_ex_aoa_imp_a} \\
		\gaoa{\vI(t)} &\definedas \Vaoa-\twonorm{\vI(t)}, \label{eq:stc_ex_aoa_imp_b} \\
    	\faoa{\vI(t)}{\qIB(t)} &\definedas \cos{\aoamax}\twonorm{\vI(t)}+\dotprod{\bvec{e}_1}{\cIB\big(\qIB(t)\big)\vI(t)},
	\end{align}
\end{subequations}
where~$\vI(t)\in\real^{3}$ is the inertial velocity,~$\qIB(t)\in\mathcal{S}^3\subset\real^{4}$ is the unit quaternion representing the transformation from~$\iframe$ to~$\bframe$, and~$\cIB(t)\definedas\cIB\big(\qIB(t)\big)$ is the direction cosine matrix associated with~$\qIB(t)$.

\subsubsection{Alternative Formulations} \label{sec2:stc_alt_forms}

We now consider a few alternative formulations that can be used in lieu of the continuous state-triggered constraint formulation presented above. Since mixed-integer methodologies were deemed not conducive to solving the powered descent guidance problem in real-time, we do not consider them further. To present these alternatives, the following discussion assumes a powered descent guidance problem that includes the state-triggered constraint formalized in~\eqref{eq:stc_ex_aoa_orig}. \\

\oursec{Alternative 1 - Naive Implementation} The first implementation we consider is a naive implementation that uses temporally-scheduled constraints. Due to its simplicity, this approach is arguably the most obvious implementation one might pursue when trying to capture the logical implication in~\eqref{eq:stc_ex_aoa_orig}. In this approach, the optimization problem is first solved without~\eqref{eq:stc_ex_aoa_orig}, and the resulting solution is used to determine the set of times~$\mathcal{T}_{\stcg}$ over which the trigger condition~$\stcg_{\aoa}(\cdot)<0$ is satisfied. The problem is then solved a second time, this time with the constraint condition~$\stcf_{\aoa}(\cdot\,,\cdot) \leq 0$ enforced for all~$t\in\mathcal{T}_{\stcg}$. The second solution is guaranteed to satisfy the constraint condition for all~$t\in\mathcal{T}_{\stcg}$. However, this solution may have time instances~$t\in\mathcal{T}_{\stcg}^{-}\subseteq\mathcal{T}_{\stcg}$ for which the trigger condition is \textit{no longer} satisfied, and where the constraint condition is needlessly enforced. This solution may also have time instances~$t\in\mathcal{T}_{\stcg}^{+}\subseteq\mathcal{T}_{\stcg}^{c}$ where the trigger condition \textit{is} now satisfied, and for which the implication in~\eqref{eq:stc_ex_aoa_orig} is violated. Thus, the inclusion of the constraint condition in the second attempt may alter the solution such that the triggering profile assumed at the onset is invalidated. Clearly, the second solution may be used to redefine the set~$\mathcal{T}_{\stcg}$, and the procedure can be repeated. However, since this approach obscures the relationship between the trigger and constraint conditions from the optimization, this behavior may persist indefinitely. Continuous state-triggered constraints avoid this issue by exposing the relationship between the trigger and constraint conditions to the optimization algorithm. \\

\oursec{Alternative 2 - Nonlinear Implementation} The second implementation we consider is one where~\eqref{eq:stc_ex_aoa_orig} is replaced with~$-\dotprod{\ex}{\vB(t)}/\twonorm{\vB(t)} \geq \cos{f_{\aoa}(t)}$, where~$f_{\aoa}(t)\definedas f_{\aoa}\big(\twonorm{\vB(t)}\big)$ is a nonlinear scalar-valued function that relates the maximum allowable angle of attack to the speed of the vehicle. Unlike the formulation used in~\eqref{eq:stc_ex_aoa_imp} which switches the constraint condition on and off in a binary fashion, this approach adjusts the constraint condition smoothly according to~$f_{\aoa}(\cdot)$. That is, instead of enforcing an angle of attack limit at high speeds and disabling it at low speeds, this implementation \textit{always} enforces an angle of attack limit that varies smoothly as a function of speed. The ability to specify a smooth transition profile to modulate the constraint condition may be advantageous in certain cases. For example, it may be desirable to leave the angle of attack unconstrained at low speeds, and to restrict it to a maximum of~$3\dg$ and~$5\dg$ at intermediate and high speeds, respectively. Using the STC framework, such a profile would either have to be conservatively approximated, or multiple cSTCs would have to be used. However, the nonlinear implementation has two disadvantages when compared to the cSTC approach. First, it is less intuitive since the relationship between the angle of attack and speed must be embedded in~$f_{\aoa}(t)$. This disadvantage can become particularly accentuated when the geometry of the constraint is multi-dimensional and difficult to visualize. Second, even if a satisfactory function~$f_{\aoa}(t)$ is obtained, it may be difficult to ensure it has sufficiently good numerical properties (e.g. optimization problems with high-order polynomials may be difficult to solve reliably). \\

\oursec{Alternative 3 - Multi-Phase Implementation} The third implementation we consider splits the problem into two temporal phases: the first phase includes an angle of attack constraint but no velocity constraint, and the second phase includes a velocity constraint but no angle of attack constraint. The terminal condition of the first phase is equated to the initial condition of the second, and both phases are solved \textit{simultaneously} as one optimal control problem. This multi-phase optimization approach ensures satisfaction of~\eqref{eq:stc_ex_aoa_orig} by construction, and is well suited for applications where the quantity and temporal ordering of the phases is known a priori. For example, this approach is ideal for optimizing the ascent trajectories of multi-staged rockets. However, since this approach presupposes the quantity and order of the phases, it may introduce unnecessary conservatism. In contrast, cSTCs are able to rearrange and introduce additional phases to achieve feasibility or improve optimality. An example of this behavior is presented in~\sref{sec:numerical_results}.


\subsection{Compound State-Triggered Constraints} \label{sec2:cmpstc}

The state-triggered constraint formulation presented in~\sref{sec2:stc} assumes scalar-valued trigger and constraint functions. This section introduces the primary contribution of this paper, namely a novel continuous formulation for \textit{compound} state-triggered constraints.

Compound state-triggered constraints are defined using \textit{vector-valued} trigger and constraint functions whose elements are composed into compound conditions using Boolean logic operations. The compound formulation is a generalization of the formulation given in~\eqref{eq:stc_formal_def}, and offers a straightforward procedure whereby simple trigger and constraint function elements can be composed into more complicated trigger and constraint conditions. We proceed as in~\sref{sec2:stc} by defining compound state-triggered constraints, presenting a continuous formulation, providing example applications, and considering alternative formulations.

\subsubsection{Formal Definition} \label{sec2:cmpstc_form_def}

The scalar state-triggered constraint defined in~\eqref{eq:stc_formal_def} may be generalized to one of the following compound state-triggered constraints
\begin{subequations} \label{eq:cmpstc_formal_def}
	\begin{align}
		& \qquad\qquad && \text{\textit{Or}-Trigger with \textit{Or}-Constraint:}  && \bigvee_{i=1}^{\stcng}{\Big(\stcg_i(\stcx) < 0\Big)}   & \;\Rightarrow\;\quad \bigvee_{i=1}^{\stcnf}{\Big(\stcf_i(\stcx) = 0\Big)}\,, && \qquad\qquad & \label{eq:cmpstc_formal_def_or_or} \\
        & \qquad\qquad && \text{\textit{Or}-Trigger with \textit{And}-Constraint:} && \bigvee_{i=1}^{\stcng}{\Big(\stcg_i(\stcx) < 0\Big)} &  \;\Rightarrow\;\quad \bigwedge_{i=1}^{\stcnf}{\Big(\stcf_i(\stcx) = 0\Big)}\,, && \qquad\qquad & \label{eq:cmpstc_formal_def_or_and} \\
        & \qquad\qquad && \text{\textit{And}-Trigger with \textit{Or}-Constraint:} && \bigwedge_{i=1}^{\stcng}{\Big(\stcg_i(\stcx) < 0\Big)} &  \;\Rightarrow\;\quad \bigvee_{i=1}^{\stcnf}{\Big(\stcf_i(\stcx) = 0\Big)}\,, && \qquad\qquad & \label{eq:cmpstc_formal_def_and_or}\\
        & \qquad\qquad && \text{\textit{And}-Trigger with \textit{And}-Constraint:} && \bigwedge_{i=1}^{\stcng}{\Big(\stcg_i(\stcx) < 0\Big)} &  \;\Rightarrow\;\quad \bigwedge_{i=1}^{\stcnf}{\Big(\stcf_i(\stcx) = 0\Big)}\,, && \qquad\qquad & \label{eq:cmpstc_formal_def_and_and}
    \end{align}
\end{subequations}
where~$\stcg(\cdot):\real^{\stcnx}\rightarrow\real^{\stcng}$ and~$\stcf(\cdot):\real^{\stcnx}\rightarrow\real^{\stcnf}$ are redefined as \textit{vector-valued} trigger and constraint functions,  and where~$\stcg_i(\stcx)$ and~$\stcf_i(\stcx)$ represent the~$\ith{i}{th}$ entries of the vectors~$\stcg(\stcx)$ and~$\stcf(\stcx)$, respectively. When \textit{at least one} element of the compound trigger condition is satisfied,~\eqref{eq:cmpstc_formal_def_or_or} implies that \textit{at least one} element of the compound constraint condition is enforced, whereas~\eqref{eq:cmpstc_formal_def_or_and} implies that \textit{all} elements of the compound constraint condition are enforced. When \textit{all} elements of the compound trigger condition are satisfied,~\eqref{eq:cmpstc_formal_def_and_or} implies that \textit{at least one} element of the compound constraint condition is enforced, whereas~\eqref{eq:cmpstc_formal_def_and_and} implies that \textit{all} elements of the compound constraint condition are enforced. Note that the compound generalizations in~\eqref{eq:cmpstc_formal_def} recover the scalar definition in~\eqref{eq:stc_formal_def} when~$\stcng=\stcnf=1$.

\subsubsection{Continuous Formulation} \label{sec2:cmpstc_cont_form}

The projected continuous formulations corresponding to the compound state-triggered constraints defined in~\eqref{eq:cmpstc_formal_def} are tabulated in Table~\ref{tab:compound_conditions}. Note that~$\stcshat_i(\stcx)\definedas-\min\big(\stcg_{i}(\stcx),0\big)\geq 0$ for all $i\in\{1,\ldots,\stcng\}$, and that~$\stch_{\land\land}(\cdot)=0$ and~$\stch_{\lor\land}(\cdot)=0$ may be implemented as standalone equality constraints \textit{only if}~$\stcf_{i}(\stcx)\geq 0$ holds for all~$i\in\{1,\ldots,\stcnf\}$. However, when element-wise non-negativity of~$\stcf(\stcx)$ does not hold,~$\stch_{\land\land}(\cdot)$ and~$\stch_{\lor\land}(\cdot)$ must be implemented as~$\stcnf$ distinct equality constraints respectively given by~$\stcf_{j}(\stcx)\cdot\prod_{i=1}^{\stcng}{\stcshat_{i}(\stcx)}=0$ and~$\stcf_{j}(\stcx)\cdot\sum_{i=1}^{\stcng}{\stcshat_{i}(\stcx)}=0$ for all~$j\in\{1,\ldots,\stcnf\}$. The standard continuous formulation can be obtained from Table~\ref{tab:compound_conditions} by replacing~$\stcshat_{i}(\stcx)$ with elements of a non-negative vector-valued~$\stcs\in\real^{\stcng}_{+}$ that adheres to~\eqref{eq:stc_cont_std_a} and~\eqref{eq:stc_cont_std_b} for each~$i\in\{1,\ldots,\stcng\}$.

\vspace{0.25cm}
\begin{table}[h!]
	\begin{center}
		\caption{Summary of Compound Continuous State-Triggered Constraints}
        \begin{tabular}{lll}
            \hhline{===}
\vspace{-0.4cm} \\
Compound  & Compound & Continuous Implementation \\
Trigger Condition & Constraint Condition  & $\stcshat_i(\stcx) \definedas  -\min\big(\stcg_i(\stcx),0\big)$ \\[0.8ex]
\hline \\
$\begin{aligned}
	\bigwedge_{i=1}^{\stcng} \Big(\stcg_i(\stcx) &< 0\Big) \\[1.5ex]
    \bigwedge_{i=1}^{\stcng} \Big(\stcg_i(\stcx) &< 0\Big) \\[1.5ex]
    \bigvee_{i=1}^  {\stcng} \Big(\stcg_i(\stcx) &< 0\Big) \\[1.5ex]
    \bigvee_{i=1}^  {\stcng} \Big(\stcg_i(\stcx) &< 0\Big)
\end{aligned}$ &
$\begin{aligned}
   \bigwedge_{i=1}^{\stcnf} \Big(\stcf_i(\stcx) &= 0\Big) \\[1.5ex]
   \bigvee_{i=1}^  {\stcnf} \Big(\stcf_i(\stcx) &= 0\Big) \\[1.5ex]
   \bigwedge_{i=1}^{\stcnf} \Big(\stcf_i(\stcx) &= 0\Big) \\[1.5ex]
   \bigvee_{i=1}^  {\stcnf} \Big(\stcf_i(\stcx) &= 0\Big)
\end{aligned}$ &
$\begin{aligned}
	&\stch_{\land \land}(\stcx) \definedas \Bigg[ \prod_{i=1}^{\stcng}{\stcshat_i(\stcx)} &\hspace{-0.3cm}&\Bigg] \cdot \Bigg[ \sum_{i=1}^{\stcnf}{\stcf_i(\stcx)}  &\hspace{-0.33cm}&\Bigg] = 0\,, \; \stcf(\stcx) \geq 0 \\[1.5ex]
    &\stch_{\land \lor}(\stcx) \definedas \Bigg[ \prod_{i=1}^{\stcng}{\stcshat_i(\stcx)} &\hspace{-0.3cm}&\Bigg] \cdot \Bigg[ \prod_{i=1}^{\stcnf}{\stcf_i(\stcx)} &\hspace{-0.33cm}&\Bigg] = 0 \\[1.5ex]
    &\stch_{\lor \land}(\stcx) \definedas \Bigg[ \sum_{i=1}^{\stcng}{\stcshat_i(\stcx)}  &\hspace{-0.3cm}&\Bigg] \cdot \Bigg[ \sum_{i=1}^{\stcnf}{\stcf_i(\stcx)}  &\hspace{-0.33cm}&\Bigg] = 0\,, \; \stcf(\stcx) \geq 0 \\[1.5ex]
    &\stch_{\lor \lor}(\stcx) \definedas \Bigg[ \sum_{i=1}^{\stcng}{\stcshat_i(\stcx)}  &\hspace{-0.3cm}&\Bigg] \cdot \Bigg[ \prod_{i=1}^{\stcnf}{\stcf_i(\stcx)} &\hspace{-0.33cm}&\Bigg] = 0
\end{aligned}$ \\
&&\vspace{-0.2cm}\\
\hhline{===}
            \label{tab:compound_conditions}
        \end{tabular}
  \end{center}
  \vspace{-1.0cm}
\end{table}

\subsubsection{Example Application} \label{sec2:cmpstc_ex}

\begin{figure}[t!]
	\begin{subfigure}{0.49\textwidth}
    	\centering
		\begin{tikzpicture}

    \tikzmath{\t  = 0.25;
              \h  = 2.00;
              \l1 = 1.00;
              \l2 = 2.50;
              \px1 = 0.5;      \py1 = 0;
              \px2 = \px1+\l1; \py2 = \py1-\t;
              \px3 = \px2-\t;  \py3 = \py2-\h+\t+\t;
              \px4 = \px3+\l2; \py4 = \py3-\t;
              \px5 = \px4+\t;  \py5 = \py4+\h;
              \px6 = \px5+\l1; \py6 = \py5-\t;
              \lx1 = \px1; \ly1 = \py1;
              \lx2 = \px2; \ly2 = \ly1;
              \lx3 = \lx2; \ly3 = \py3;
              \lx4 = \px4; \ly4 = \ly3;
              \lx5 = \lx4; \ly5 = \py5;
              \lx6 = \px6; \ly6 = \ly5;
              \cx = 0.5*(\lx3+\lx4); \cy = \ly3;
              \tx = \lx4+0.0*\l1; \ty=\ly4-0.5*\h;
              \dy1 = \py1+0.1;
              \dy2 = \dy1+0.4;
              \dy3 = 0.5*(\dy1+\dy2);
              \dxl3 = 0.5*(\lx2+\cx);
              \dxr3 = 0.5*(\cx+\lx5);
              \dx4 = \lx6+0.1;
              \dx5 = \dx4+0.4;
              \dx6 = 0.5*(\dx4+\dx5);
              \dy6 = 0.5*(\ly4+\ly5);}

    \fill[pattern=north west lines, pattern color=black] (\px1,\py1) rectangle (\px2,\py2);
    \fill[pattern=north west lines, pattern color=black] (\px2,\py2) rectangle (\px3,\py3);
    \fill[pattern=north west lines, pattern color=black] (\px3,\py3) rectangle (\px4,\py4);
    \fill[pattern=north west lines, pattern color=black] (\px4,\py4) rectangle (\px5,\py5);
    \fill[pattern=north west lines, pattern color=black] (\px5,\py5) rectangle (\px6,\py6);
    
    \draw[black,cap=round] (\lx1,\ly1) -- (\lx2,\ly2);
    \draw[black,cap=round] (\lx2,\ly2) -- (\lx3,\ly3);
    \draw[black,cap=round] (\lx3,\ly3) -- (\lx4,\ly4);
    \draw[black,cap=round] (\lx4,\ly4) -- (\lx5,\ly5);
    \draw[black,cap=round] (\lx5,\ly5) -- (\lx6,\ly6);
    
    \draw[black] (\lx2,\dy1) -- (\lx2,\dy2);
    \draw[black] (\lx5,\dy1) -- (\lx5,\dy2);
    \draw[black] (\cx,\dy1) -- (\cx,\dy2);
    \draw[black,<->] (\lx2,\dy3) -- (\cx,\dy3);
    \draw[black,<->] (\cx,\dy3) -- (\lx5,\dy3);
    \draw (\dxl3,\dy3) node[anchor=south] {$a$};
    \draw (\dxr3,\dy3) node[anchor=south] {$b$};
    
    \draw[black] (\dx4,\ly4) -- (\dx5,\ly4);
    \draw[black] (\dx4,\ly6) -- (\dx5,\ly6);
    \draw[black,<->] (\dx6,\ly4) -- (\dx6,\ly6);
    \draw (\dx6,\dy6) node[anchor=west] {$c$};
    
    \filldraw[black] (\cx,\cy) circle (2pt);
    \draw (\cx,\cy) to[out=-90,in=-180] (\tx,\ty);
    \draw (\tx,\ty) node[anchor=west] {Landing site};
    
    \draw (\cx-0.1,\cy+0.3) node[anchor=east] {};
    \draw[black,thick,->] (\cx,\cy) -- +(0.75,0);
    \draw[black,thick,->] (\cx,\cy) -- +(0,0.75);
    \draw (\cx+0.1,\cy2+0.7) node[anchor=west] {$\bvec{e}_1$};
    \draw (\cx+0.4,\cy2+0.2) node[anchor=west] {$\bvec{e}_2$};
\end{tikzpicture}
    	\caption{Landing in a Crater or Canyon}
    	\label{fig:canyon}
   	\end{subfigure}
    \begin{subfigure}{0.49\textwidth}
    	\centering
		\begin{tikzpicture}

    \tikzmath{\t  = 0.25;
              \mx = 1.00;
              \h1 = 2.00; \w1 = 4.00;
              \h2 = 1.00; \w2 = 1.50;
              \h3 = 0.75; \w3 = \w1-\w2;
              \px1 = 0;            \py1 = 0;
              \px2 = \px1+\mx*\h1; \py2 = \py1-\h1;
              \px3 = \px2+\w1;     \py3 = \py2;
              \px4 = \px1+\w2;     \py4 = \py1;
              \px5 = \px4+\mx*\h2; \py5 = \py4-\h2;
              \px6 = \px5+\w3;     \py6 = \py5;
              \px7 = \px4;         \py7 = \py4+\h3;
              \px8 = \px5;         \py8 = \py5+\h3;
              \px9 = \px6;         \py9 = \py6+\h3;
              \pxc = \px7+\w3;     \pyc = \py7;
              \tx1 = \px2+1;       \ty1 = \py2-0.5;
              \tx2 = \pxc+0.5;     \ty2 = \pyc-0.2;}
    
    \fill[color=black!10] (\px4,\py4) -- (\px5,\py5) -- (\px8,\py8) -- (\px7,\py7) -- cycle;
    \fill[color=black!15] (\px5,\py5) -- (\px6,\py6) -- (\px9,\py9) -- (\px8,\py8) -- cycle;
    \fill[color=black!5 ] (\px7,\py7) -- (\px8,\py8) -- (\px9,\py9) -- (\pxc,\pyc) -- cycle;
    
    \draw[color=black!100,cap=round] (\px4,\py4) -- (\px5,\py5);
    \draw[color=black!100,cap=round] (\px5,\py5) -- (\px8,\py8);
    \draw[color=black!100,cap=round] (\px8,\py8) -- (\px7,\py7);
    \draw[color=black!100,cap=round] (\px4,\py4) -- (\px7,\py7);
    
    \draw[color=black!100,cap=round] (\px5,\py5) -- (\px6,\py6);
    \draw[color=black!100,cap=round] (\px6,\py6) -- (\px9,\py9);
    \draw[color=black!100,cap=round] (\px9,\py9) -- (\px8,\py8);
    
    \draw[color=black!100,cap=round] (\pxc,\pyc) -- (\px9,\py9);
    \draw[color=black!100,cap=round] (\pxc,\pyc) -- (\px7,\py7);

    \draw[black,cap=round] (\px4,\py4) -- (\px1,\py1);
    \draw[black,cap=round] (\px1,\py1) -- (\px2,\py2);
    \draw[black,cap=round] (\px2,\py2) -- (\px3,\py3);
    \draw[black,cap=round] (\px3,\py3) -- (\px6,\py6);
    
    \draw[black] (\px1-\mx*0.1,\py1+0.1) -- (\px1-\mx*0.3,\py1+0.3);
	\draw[black] (\px4-\mx*0.1,\py4+0.1) -- (\px4-\mx*0.3,\py4+0.3);
    \draw[black,<->] (\px1-\mx*0.2,\py1+0.2) -- (\px4-\mx*0.2,\py4+0.2);
    
    \draw[black] (\px3+0.1,\py3) -- (\px3+0.5,\py3);
    \draw[black] (\px6+0.1,\py6) -- (\px6+0.5,\py6);
    \draw[black] (\px9+0.1,\py9) -- (\px9+0.5,\py9);
    \draw[black,<->] (\px3+0.3,\py3) -- (\px6+0.3,\py6);
    \draw[black,<->] (\px6+0.3,\py6) -- (\px9+0.3,\py9);
    
    \draw[black,thick,->] (\px2,\py2) -- (\px2+0.75,\py2);
    \draw[black,thick,->] (\px2,\py2) -- (\px2-\mx*0.5,\py2+0.5);
    \draw[black,thick,->] (\px2,\py2) -- (\px2,\py2+0.75);
    
    \filldraw[black] (\px2,\py2) circle (2pt);
    \draw (\px2-0.1,\py2-0.2) node[anchor=east] {};
    \draw (\px2,\py2) to[out=-90,in=-180] (\tx1,\ty1);
    \draw (\tx1,\ty1) node[anchor=west] {Landing site};
    
    \draw (\px7*0.25+\px8*0.25+\px9*0.25+\pxc*0.25,
           \py7*0.25+\py8*0.25+\py9*0.25+\pyc*0.25) to[out=45,in=-180] (\tx2,\ty2);
    \draw (\tx2,\ty2) node[anchor=west] {Keep-out volume};
    
    \draw (\px2+0.1,\py2+0.6) node[anchor=west] {$\bvec{e}_1$};
    \draw (\px2+0.5,\py2+0.2) node[anchor=west] {$\bvec{e}_2$};
    \draw (\px2-0.4,\py2+0.3) node[anchor=east] {$\bvec{e}_3$};
    
    \draw (\px1*0.5+\px4*0.5-\mx*0.3,\py1*0.5+\py4*0.5+0.3) node[anchor=south] {$a$};
    \draw (\px3*0.5+\px6*0.5+0.5,\py3*0.5+\py6*0.5) node[anchor=west] {$b$};
    \draw (\px6*0.5+\px9*0.5+0.5,\py6*0.5+\py9*0.5) node[anchor=west] {$c$};
\end{tikzpicture}
    	\caption{Landing Near Cliffs or Mesas}
    	\label{fig:mesa}
   	\end{subfigure}
    \caption{An example powered descent guidance application where compound state-triggered constraints  are used to avoid collisions with geological formations near the landing site.}
   	\label{fig:cmpstc_ex}
   	\vspace{-0.5cm}
\end{figure}
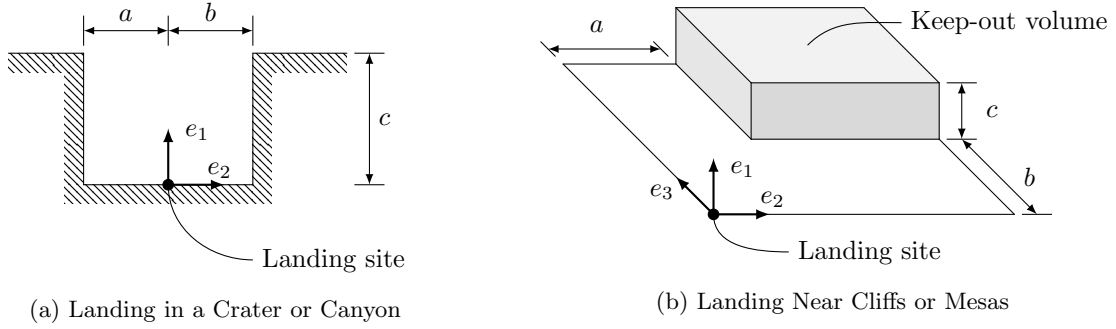

To motivate the use of compound STCs, consider a scenario where the objective is to land in the vicinity of large geologic formations such as the craters, canyons, cliffs, or mesas. Figure~\ref{fig:canyon} illustrates a scenario where the landing site is located at the bottom of a crater or canyon. Using the parameters defined in the figure, the collision avoidance constraint is given by the following \textit{or}-triggered compound STC
\begin{equation*}
	\big(\dotprod{\bvec{e}_2}{\rI(t)} < -a \big) \lor \big(b < \dotprod{\bvec{e}_2}{\rI(t)}\big) \;\Rightarrow\; c-\dotprod{\bvec{e}_1}{\rI(t)} \leq 0,
\end{equation*}
where~$\rI(t)\in\real^{3}$ is the vehicle's inertial position vector. The corresponding compound cSTC is given by
\begin{subequations}
	\begin{align*}
		\hor{\rI(t)} &\definedas \bigg[ \sum_{i=1}^{2}{-\min\Big(\gor{i}{\rI(t)} , 0 \Big)} \bigg] \cdot \for{\rI(t)} \leq 0, \\
		\gor{1}{\rI(t)} &\definedas a + \dotprod{\bvec{e}_2}{\rI(t)}, \\
    	\gor{2}{\rI(t)} &\definedas b - \dotprod{\bvec{e}_2}{\rI(t)}, \\
    	\for{\rI(t)}    &\definedas c - \dotprod{\bvec{e}_1}{\rI(t)}.
	\end{align*}
\end{subequations}

Figure~\ref{fig:mesa} illustrates a similar but more challenging scenario, where a keep-out volume is used to steer the vehicle away from geological formations such as cliffs or mesas. Following the procedure outlined above, we obtain a continuous formulation of the collision avoidance constraints using the following \textit{and}-triggered compound cSTC
\begin{subequations} \label{eq:cmpstc_ex_and}
	\begin{align}
		\hand{\rI(t)} &\definedas \bigg[ \prod_{i=1}^{2}{-\min\Big(\gand{i}{\rI(t)} , 0 \Big)} \bigg] \cdot \fand{\rI(t)} \leq 0, \\
		\gand{1}{\rI(t)} &\definedas a - \dotprod{\bvec{e}_2}{\rI(t)}, \\
    	\gand{2}{\rI(t)} &\definedas b - \dotprod{\bvec{e}_3}{\rI(t)}, \\
    	\fand{\rI(t)} &\definedas c - \dotprod{\bvec{e}_1}{\rI(t)}.
	\end{align}
\end{subequations}

\subsubsection{Alternative Formulations} \label{sec2:cmpstc_alt_form}

The scenario described in~\sref{sec2:cmpstc_ex} can instead be addressed using the glide-slope cone constraint formulated in~\cite{behcetjgcd07}. However, enforcing a sufficiently steep constraint over the entire trajectory may be overly conservative, and may prohibitively restrict the set of feasible initial conditions (e.g., Mars' Valles Marineris reaches depths of up to 7 km). Alternatively, enforcing the constraint over only a terminal time interval may allow the trajectory to violate the physical constraints at times prior to said interval. While this approach may work in certain cases, the heuristics it employs may fall short for applications that have a large set of possible initial conditions.

Alternatively, this scenario can be addressed using the multi-phase optimization approach presented in~\cite{bhasin2016fuel}, which is similar to \textit{Alternative 3} discussed in~\sref{sec2:stc_alt_forms}. The multi-phase optimization approach is well suited for this example application since the quantity and ordering of the phases is known a priori (provided that the altitude decreases monotonically). Nevertheless, we argue that compound state-triggered constraints offer a simple way to formulate complex keep-out volumes in the scenario described above.


\newpage
\section{Problem Statement \& Convex Formulation} \label{sec:prob_state}

In this section, we outline a problem statement for a non-convex 6-DoF powered descent guidance problem, and provide a high-level description of the successive convexification procedure used to solve it. The content of this section is drawn from Sections II and III of~\cite{SzmukReynolds2018}, and is provided as context for~\sref{sec:numerical_results}.

\subsection{Non-Convex Problem Statement} \label{sec3:ncvx}

We consider the minimum-time 6-DoF powered descent guidance problem given in Problem~\ref{prob:ncvx}, which is comprised of nonlinear dynamics and non-convex state and control constraints. This problem formulation includes the baseline problem formulation, free-ignition-time modification, and ellipsoidal aerodynamics model detailed in Sections II.A-II.E of~\cite{SzmukReynolds2018}. Additionally, we include the velocity-triggered angle of attack constraint discussed in~\sref{sec2:stc_ex}, and the collision avoidance constraint formulated in~\sref{sec2:cmpstc_ex}.

Our problem formulation consists of three time epochs: \textit{initial time}~$\tin$, \textit{ignition time}~$\tig$, and \textit{final time}~$\tf$. These epochs are defined such that~$\tin\leq\tig<\tf$. We call the time interval~$t\in\intee{\tin}{\tig}$ the \textit{coast phase}, and the interval~$t\in\intee{\tig}{\tf}$ the \textit{burn phase}. During the coast phase, the engine is off, and the vehicle state evolves passively along a predetermined trajectory (as is discussed shortly). We denote the \textit{coast time} by~$\tc\definedas\tig-\tin$, and the maximum allowable coast time by~$\tcmax$. During the burn phase, the engine is on, and the vehicle actively maneuvers to achieve its landing objective. We denote the \textit{burn time} by~$\tb\definedas\tf-\tig$.

The states of the problem are the mass~$m(t)\in\real_{++}$, angular velocity~$\omegaB(t)\in\real^{3}$, and the position~$\rI(t)$, velocity~$\vI(t)$, and unit quaternion~$\qIB(t)$ introduced in~\sref{sec2:stc_ex} and~\sref{sec2:cmpstc_ex}. The control of the problem is the thrust vector~$\TB(t)\in\real^{3}$. These are summarized as the state vector
$$\xx(t)\definedas\big[\m(t)\;\rI^T(t)\;\vI^T(t)\;\qIB^T(t)\;\omegaB^T(t)\big]^T\in\real^{14},$$
and the control vector~$\uu(t)\definedas\TB(t)\in\real^{3}$. Since this is a free-final-time minimum-time problem with a free-ignition-time modification, the cost function of the problem is given by $J(\zz)\definedas\tb$, and the solution variable is given by
\begin{equation*}
    \zz(t)\definedas\big[\tc\;\tb\;\xx^T(t)\;\uu^T(t)\big]^T\in\real^{19}.
\end{equation*}

At the ignition time epoch, the following boundary conditions are enforced
\begin{equation*} \label{eq:ics}
    \m(\tig) = \mig,\quad \rI(\tig)=\prig(\tc),\quad \vI(\tig)=\pvig(\tc),\quad \omegaB(\tig)=\bvec{0},
\end{equation*}
where, due the free-ignition-time modification, the ignition time position and velocity are constrained to lie on curves described by the polynomials~$\prig:\real\rightarrow\real^{3}$ and~$\pvig:\real\rightarrow\real^{3}$. Although our problem formulation considers aerodynamic effects, we assume that~$\prig(\cdot)$ and~$\pvig(\cdot)$ describe an aerodynamics-free free fall trajectory given by
\begin{equation*} \label{eq:fit}
    \prig(\tc) \definedas \rIin + \vIin\,\tc + \frac{1}{2}\gI\tc^2, \quad \pvig(\tc) \definedas \vIin + \gI\tc,
\end{equation*}
where~$\rIin$ and~$\vIin$ are prescribed initial-time position and velocity vectors, and~$\gI\in\real^{3}$ is the (constant) gravity vector. We note, however, that~$\prig(\cdot)$ and~$\pvig(\cdot)$ can be formulated to reflect higher order effects such as aerodynamics. At the final time epoch, the following boundary conditions are enforced
\begin{equation*} \label{eq:tcs}
    \rI(\tf) = \bvec{0},\quad \vI(\tf)=\bvec{0},\quad \qIB(\tf)=\qidentity,\quad \omegaB(\tf)=\bvec{0},
\end{equation*}
where~$\qidentity\definedas[1\;\,0\;\,0\;\,0]^T$ denotes the identity quaternion.

The mass state of the vehicle evolves according to the affine (in thrust) mass depletion dynamics used in~\cite{szmuk2016successive}, while the translational and attitude states evolve according to 6-DoF rigid body dynamics. The former depends on the mass depletion coefficients~$\mdotalpha\in\real_{++}$ and~$\mdotbeta\in\real_{+}$ that are functions of the specific impulse, nozzle exit area, and atmospheric back pressure. The latter is affected by the commanded thrust~$\TB(t)$ and the aerodynamic force~$\aeroB(t)\in\real^{3}$, which apply torque to the vehicle with moment arms~$\rTB\in\real^{3}$ and $\rCPB\in\real^{3}$, respectively. The aerodynamic force is defined by the ellipsoidal aerodynamic model introduced in Section II.E of~\cite{SzmukReynolds2018}. Additionally, we use~$\Omega(\cdot)$ to denote the $4\times 4$ skew-symmetric matrix associated with the quaternion kinematics, and~$\inertia\in\spd{3}$ to denote the vehicle's body-fixed constant inertia.

The problem formulation includes state inequality constraints that keep the mass above a dry mass $\mdry\in\real_{++}$, the position inside a glide slope cone with a half-angle of~$\glideslope\in\intei{0\dg}{90\dg}$, the tilt angle less than $\tiltmax\in\intie{0\dg}{90\dg}$, and the angular velocity magnitude less than~$\omegamax\in\real_{++}$. The matrices in the glide slope and tilt constraints are defined as~$\Hgs\definedas[\bvec{e}_2\;\,\bvec{e}_3]^T\in\real^{2\times 3}$ and~$\Htilt\definedas[\bvec{e}_3\;\,\bvec{e}_4]^T\in\real^{2\times 4}$, respectively. The formulation also includes two control constraints. The first constraint limits the thrust magnitude to a minimum value of~$\Tmin$ and a maximum value of~$\Tmax$, while the second ensures that the gimbal angle does not exceed~$\gimbalmax\in\intii{0\dg}{90\dg}$. We refer the reader to~\cite{SzmukReynolds2018} for more details.

\boxing{t!}{problem}{prob:ncvx}{16cm}{Non-Convex Optimal Control Problem}{
 	\begin{tabular}{llll}
 		\underline{Cost Function}: &&& \\
 		&& $ {\begin{aligned}
        		&\underset{\tc,\,\tb,\,\TB(t)}{\text{minimize}}\;\; \tb \\[1ex]
                &\text{s.t.}\quad\tc\in\intee{0}{\tcmax}
            \end{aligned}}$
        & \\
 		\underline{Boundary Conditions}: &&& \\
 		&& $ {\begin{aligned}
 		                 \m(\tig) &= \mig       && & \qIB(\tf)     &= \qIBf \\ 
 				        \rI(\tig) &= \prig(\tc) && & \rI(\tf)      &= \rIf \\
 					    \vI(\tig) &= \pvig(\tc) && & \vI(\tf)      &= \vIf \\
 					\omegaB(\tig) &= \omegaBig  && & \omegaB(\tf)  &= \omegaBf
         	\end{aligned}} $ 
         & \\
 		\underline{Dynamics}: &&& \\
 		&& $ {\begin{aligned}
 					\dot{m}(t) &= -\mdotalpha\twonorm{\TB(t)} - \mdotbeta \\
 					\rIdot(t) &= \vI(t) \\
 					\vIdot(t) &= \frac{1}{m(t)} \cBI(t) \big(\TB(t) + \aeroB(t)\big) + \gI \\
 					\qIBdot(t) &= \frac{1}{2}\OMEGA{\omegaB(t)}\qIB(t) \\
                    \inertia \omegaBdot(t) &= \crossprod{\rTB}{\TB(t)} + \crossprod{\rCPB}{\aeroB(t)}- \crossprod{\omegaB(t)}{\inertia\omegaB(t)}
 				\end{aligned}} $ 
         & \\
        $ {\begin{aligned}
        		&\text{\underline{State Constraints}:} \\
                & \\
                & \\
                & \\
                & \\
		 		&\text{\underline{Control Constraints}:} \\
                & \\
                & \\
 				&\text{\underline{State-Triggered Constraints}:} \\
				& \\
             \end{aligned}} $
        && $ {\begin{aligned}
         		& \\
                \mdry &\leq m(t) \\
 				\tan{\glideslope}\twonorm{\Hgs\rI(t)} &\leq \dotprod{\ex}{\rI(t)} \\
 				\cos{\tiltmax} &\leq 1-2\twonorm{\Htilt\qIB(t)} \\
 				\twonorm{\omegaB(t)} &\leq \omegamax \\
 				& \\
                0 < \Tmin \leq \twonorm{\TB(t)} &\leq \Tmax \\
				\cos{\gimbalmax}\twonorm{\TB(t)} &\leq \dotprod{\ex}{\TB(t)} \\
                & \\
	         	\haoa{\vI(t)}{\qIB(t)} &\leq 0 \\
	         	\hand{\rI(t)} &\leq 0
             \end{aligned}} $
        & $ {\begin{aligned}
        		& \\
        		& \\
        		& \\
        		& \\
        		& \\
        		& \\
         		& \\
         		& \\
         		& \\
                & \hspace{-0.5cm}\text{See}\; \eqref{eq:stc_ex_aoa_imp} \\
                & \hspace{-0.5cm}\text{See}\; \eqref{eq:cmpstc_ex_and}
             \end{aligned}} $ \\
         \null
 	\end{tabular}
}


\subsection{Convex Formulation} \label{sec3:cvx}

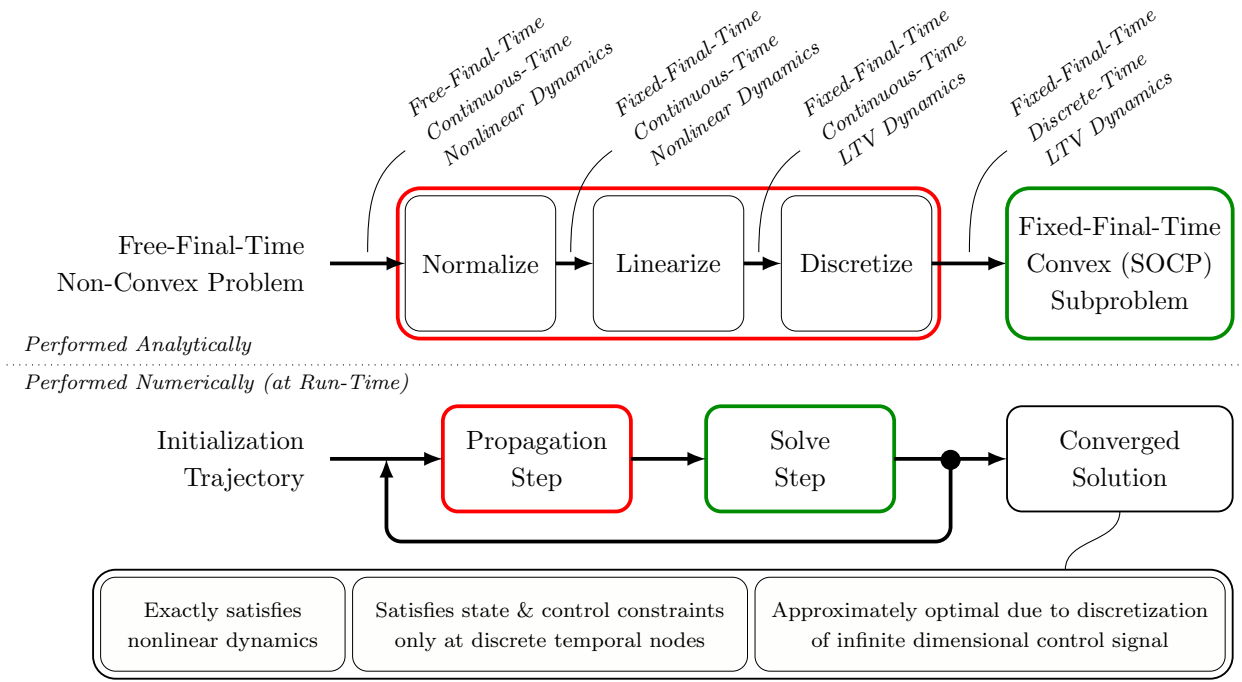
\begin{figure}[t!]
	\centering
  	\begin{tikzpicture}
    \definecolor{propcolor}{rgb}{1.0,0.0,0.0}
    \definecolor{solvecolor}{rgb}{0.0,0.55,0.0}
    \definecolor{fillcolor123}{rgb}{1,1,1}
    \colorlet{filltextcolor}{beige!10};
    \colorlet{linecolor123}{black!100};
    \colorlet{fillcolor}{beige!0};
    \colorlet{arrowcolor}{black!100};

    \newcommand{\lwarrow}{0.50mm};
    \newcommand{\lwthin}{0.10mm};
    \newcommand{\lwthick}{0.50mm};
    \newcommand{\lwtext}{0.25mm};
    \newcommand{\textrot}{35};

    \tikzmath{\cx1=0.00; \cy1=0;    \w1=2.0; \h1=1.8; \c1=0.2;
              \cx2=2.50; \cy2=\cy1; \w2=\w1; \h2=\h1; \c2=\c1;
              \cx3=5.00; \cy3=\cy1; \w3=\w1; \h3=\h1; \c3=\c1;
              \cx4=8.50; \cy4=\cy1; \w4=3.0; \h4=\h1; \c4=\c1;
              \cx0=-2.0; \cy0=\cy1; \w0=3.5; \h0=\h1; \c0=\c1;
              \dx1=0.1; \dy1=0.1; \dc1=\c1+0.1;
              \ccx0=0.5*\cx0+0.5*\cx1-0.25*\w1;
              \ccx1=0.5*\cx1+0.25*\w1+0.5*\cx2-0.25*\w2;
              \ccx2=0.5*\cx2+0.25*\w2+0.5*\cx3-0.25*\w3;
              \ccx3=0.5*\cx3+0.25*\w3+0.5*\cx4-0.25*\w4;
              \ddx0=-0.05; \ddx1=0.50;
              \ddy0= 0.15; \ddy1=1.50;
              \cx5=\cx1+0.75; \cy5=\cy1-2.6; \w5=\w4-0.5; \h5=1.4; \c5=\c4;
              \cx6=\cx1+4.25; \cy6=\cy5;     \w6=\w5;     \h6=\h5; \c6=\c5;
              \cxt=\cx1-3.5; \cyt=0.5*\cy1+0.5*\cy5-0.05; \dcyt=0.25;
              \cx7=\cx4; \cy7=\cy5; \w7=\w4; \h7=\h5; \c7=\c5;
              \cx8=\cx0; \cy8=\cy5; \w8=\w0; \h8=\h5; \c8=\c5;
              \ppx1=0.5*\cx6+0.25*\w6+0.5*\cx7-0.25*\w7; \ppy1=\cy6;
              \ppx2=\ppx1;                               \ppy2=\ppy1-1.1;
              \ppx3=0.5*\cx5-0.25*\w5+0.5*\cx8;          \ppy3=\ppy2;
              \ppx4=\ppx3;                               \ppy4=\ppy1;
              \ppc=0.2;
              \ww1=3.25; \ww2=5.25; \ww3=6.25; \hh=1.25; \ddw=0.1; \ccw=\c5; \ddhh=-2.2;
              \ttcx1=\cx7+0.5*\w7-\ww3-\ww2-0.5*\ww1-3*\ddw;
              \ttcx2=\cx7+0.5*\w7-\ww3-0.5*\ww2-2*\ddw;
              \ttcx3=\cx7+0.5*\w7-0.5*\ww3-\ddw;
              \ttcy=\cy7+\ddhh;
      }
    
    \draw (\cx0,\cy0) node[anchor=east] {\begin{tabular}{r}
                                             Free-Final-Time \\
                                             Non-Convex Problem
                                         \end{tabular}};
    
    \cpane{\cx1-0.5*\w1-\dx1}{\cx3+0.5*\w3+\dx1}
          {\cy1+0.5*\h1+\dy1}{\cy3-0.5*\h3-\dy1}{0}{\dc1}{color=propcolor,line width=\lwthick,fill=fillcolor};

    \draw[arrowcolor,line width=\lwarrow,->] (\cx0,\cy0) -- (\cx1-0.5*\w1,\cy1);

    \pane{\cx1}{\cy1}{0}{\w1}{\h1}{\c1}{color=linecolor123,line width=\lwthin,fill=fillcolor123};
    \draw (\cx1,\cy1) node[anchor=center] {Normalize};
    
    \draw[arrowcolor,line width=\lwarrow,->] (\cx1+0.5*\w1,\cy1) -- (\cx2-0.5*\w2,\cy2);
    
    \pane{\cx2}{\cy2}{0}{\w2}{\h2}{\c2}{color=linecolor123,line width=\lwthin,fill=fillcolor123};
    \draw (\cx2,\cy2) node[anchor=center] {Linearize};
    
    \draw[arrowcolor,line width=\lwarrow,->] (\cx2+0.5*\w2,\cy2) -- (\cx3-0.5*\w3,\cy3);

    \pane{\cx3}{\cy3}{0}{\w3}{\h3}{\c3}{color=linecolor123,line width=\lwthin,fill=fillcolor123};
    \draw (\cx3,\cy3) node[anchor=center] {Discretize};
    
    \draw[arrowcolor,line width=\lwarrow,->] (\cx3+0.5*\w3,\cy3) -- (\cx4-0.5*\w4,\cy4);

    \pane{\cx4}{\cy4}{0}{\w4}{\h4+2*\dy1}{\dc1}{color=solvecolor,line width=\lwthick,fill=fillcolor};
    \draw (\cx4,\cy4) node[anchor=center] {\begin{tabular}{c}
                                             Fixed-Final-Time \\
                                             Convex (SOCP) \\
                                             Subproblem
                                           \end{tabular}};
    
    \draw[black] (\ccx0+0*\ddx0,\cy0+\ddy0) to[out=90,in=180+\textrot] (\ccx0+\ddx1,\cy0+\ddy1) node[anchor=west,rotate=\textrot]{\it\footnotesize
        \begin{tabular}{l}
            Free-Final-Time \\
            Continuous-Time \\
            Nonlinear Dynamics
        \end{tabular}};
	
	\draw[black] (\ccx1+\ddx0,\cy1+\ddy0) to[out=90,in=180+\textrot] (\ccx1+\ddx1,\cy1+\ddy1) node[anchor=west,rotate=\textrot]{\it\footnotesize
        \begin{tabular}{l}
            Fixed-Final-Time \\
            Continuous-Time \\
            Nonlinear Dynamics
        \end{tabular}};
    
    \draw[black] (\ccx2+\ddx0,\cy2+\ddy0) to[out=90,in=180+\textrot] (\ccx2+\ddx1,\cy2+\ddy1) node[anchor=west,rotate=\textrot]{\it\footnotesize
        \begin{tabular}{l}
            Fixed-Final-Time \\
            Continuous-Time \\
            LTV Dynamics
        \end{tabular}};
    
    \draw[black] (\ccx3+0*\ddx0,\cy3+\ddy0) to[out=90,in=180+\textrot] (\ccx3+\ddx1,\cy3+\ddy1) node[anchor=west,rotate=\textrot]{\it\footnotesize
        \begin{tabular}{l}
            Fixed-Final-Time \\
            Discrete-Time \\
            LTV Dynamics
        \end{tabular}};
    
    \draw (\cxt-2.7,\cyt+\dcyt) node[anchor=west] {\it\footnotesize Performed Analytically};
    \draw[black,dotted] (\cx1-6.3,\cyt) --+(16.5,0.0);
    \draw (\cxt-2.7,\cyt-\dcyt) node[anchor=west] {\it\footnotesize Performed Numerically (at Run-Time)};
    
    \pane{\cx5}{\cy5}{0}{\w5}{\h5}{\c5}{color=propcolor,line width=\lwthick,fill=fillcolor};
    \draw (\cx5,\cy5) node[anchor=center] {\begin{tabular}{c}
                                             Propagation \\ Step
                                           \end{tabular}};
    
    \pane{\cx6}{\cy6}{0}{\w6}{\h6}{\c6}{color=solvecolor,line width=\lwthick,fill=fillcolor};
    \draw (\cx6,\cy6) node[anchor=center] {\begin{tabular}{c}
                                             Solve \\ Step
                                           \end{tabular}};
    
    \pane{\cx7}{\cy7}{0}{\w7}{\h7}{\c7}{color=black,line width=\lwtext,fill=fillcolor};
    \draw (\cx7,\cy7) node[anchor=center] {\begin{tabular}{c}
                                             Converged \\
                                             Solution
                                           \end{tabular}};
    
    \draw (\cx8,\cy8) node[anchor=east] {\begin{tabular}{r}
                                             Initialization \\
                                             Trajectory
                                         \end{tabular}};
    
    \draw[arrowcolor,line width=\lwarrow,->] (\cx5+0.5*\w5,\cy5) -- (\cx6-0.5*\w6,\cy6);
    \draw[arrowcolor,line width=\lwarrow,->] (\cx6+0.5*\w6,\cy6) -- (\cx7-0.5*\w7,\cy7);
    \draw[arrowcolor,line width=\lwarrow,->] (\cx8,\cy8) -- (\cx5-0.5*\w5,\cy5);
    \draw[arrowcolor,line width=\lwarrow,*->]
        (\ppx1,\ppy1+0.12) --
        (\ppx2,\ppy2+\ppc) to[out=-90,in=0] (\ppx2-\ppc,\ppy2) --
        (\ppx3+\ppc,\ppy3) to[out=-180,in=-90] (\ppx3,\ppy3+\ppc) --
        (\ppx4,\ppy4);
    
    \cpane{\ttcx1-0.5*\ww1-\ddw}{\ttcx3+0.5*\ww3+\ddw}
          {\ttcy+0.5*\hh+\ddw}{\ttcy-0.5*\hh-\ddw}{0}{\ccw+\ddw}{color=black,line width=\lwtext,fill=fillcolor};
    \pane{\ttcx1}{\ttcy}{0}{\ww1}{\hh}{\ccw}{color=black,line width=\lwthin,fill=filltextcolor};
    \pane{\ttcx2}{\ttcy}{0}{\ww2}{\hh}{\ccw}{color=black,line width=\lwthin,fill=filltextcolor};
    \pane{\ttcx3}{\ttcy}{0}{\ww3}{\hh}{\ccw}{color=black,line width=\lwthin,fill=filltextcolor};
    
    \draw (\ttcx1,\ttcy) node[anchor=center] {\footnotesize\begin{tabular}{c}
                                             Exactly satisfies \\
                                             nonlinear dynamics
                                         \end{tabular}};
    \draw (\ttcx2,\ttcy) node[anchor=center] {\footnotesize\begin{tabular}{c}
                                             Satisfies state \& control constraints \\
                                             only at discrete temporal nodes
                                         \end{tabular}};
    \draw (\ttcx3,\ttcy) node[anchor=center] {\footnotesize\begin{tabular}{c}
                                             Approximately optimal due to discretization \\ of infinite dimensional control signal
                                         \end{tabular}};
    \draw (\ttcx3+1.0,\ttcy+0.5*\hh+\ddw) to[out=90,in=-90] (\cx7,\cy7-0.5*\h7);

\end{tikzpicture}
  	\caption{High-level description of the successive convexification algorithm. The top half of the figure shows the three analytical steps that comprise the \textit{propagation step}. The propagation step generates a fixed-final-time second-order cone programming subproblem that is solved during the \textit{solve step}. The propagation and solve steps are executed numerically at run-time. The converged solution exhibits the three properties detailed in the text boxes at the bottom.}
  	\label{fig:convexification}
  	\vspace{-0.5cm}
\end{figure}

Successive convexification is an iterative solution methodology wherein the original non-convex problem is converted into a sequence of manageable convex approximations, or subproblems. These convex subproblems can be solved using fast and reliable interior point method algorithms~\cite{Domahidi2013ecos,dueri2016customized}, thereby bolstering the real-time capabilities of the successive convexification methodology. While the field of solving convex optimization problems is quite mature, methods that convert a non-convex optimal control problem into a numerically tractable parameter optimization problem are more varied.

Section III in~\cite{SzmukReynolds2018} outlines a procedure whereby the original free-final-time nonlinear continuous-time discrete-time optimal control problem (i.e. Problem~\ref{prob:ncvx}) is converted into a fixed-final-time linear-time-varying optimal control problem. This procedure consists of three analytical steps that \textit{normalize}, \textit{linearize}, and \textit{discretize} the problem. The result of these three steps is a numerically tractable second-order cone programming subproblem that locally approximates Problem~\ref{prob:ncvx}. Due to the linearization step, this local approximation is \textit{guaranteed} to be convex. This procedure is illustrated in Figure~\ref{fig:convexification}.

In practice, successive convexification works as follows. The iterative process is initialized with a user-specified reference trajectory. During the first iteration, the three analytical steps mentioned above are executed as part of the \textit{propagation step}. The propagation step computes a set of matrices that characterize the approximated discrete-time linear-time-varying dynamics and linearized state and control constraints associated with the convex subproblem. The resulting second-order cone programming subproblem is subsequently solved in the \textit{solve step}, and its solution is used as the reference trajectory for the subsequent successive convexification iteration.

Since the subproblem is an approximation of the non-convex elements of Problem~\ref{prob:ncvx} obtained through linearization, one might expect the process described above to generate approximate solutions. While this is true in some sense, we emphasize that a converged solution of the process shown in Figure~\ref{fig:convexification} has the following three properties: (i) the solution exactly satisfies the nonlinear dynamics, (ii) the solution satisfies the state and control constraints at a finite number of temporal nodes defined during the discretization step, and (iii) the solution is approximately (locally) optimal. The third property is due to the first-order-hold interpolation used in~\cite{SzmukReynolds2018} to represent the infinite dimensional control signal as a finite-dimensional linearly interpolated signal. The reader is referred to Section III of~\cite{SzmukReynolds2018} for more details on the propagation step and on the successive convexification algorithm.

Convergence guarantees for a successive convexification algorithm were presented in~\cite{SCvx_cdc16,SCvx_2018arXiv}. The algorithm presented in these works employed a virtual control term, hard trust regions, and an exact penalty method, and was guaranteed to converge globally to a (not necessarily feasible) stationary point. The result stated that a converged solution with zero virtual control was a local optimum of the original problem. In contrast, the algorithm used herein and in~\cite{SzmukReynolds2018} uses soft trust regions that are augmented to the cost, and does not yet have convergence results. However, it is similar in structure to the algorithm presented  in~\cite{SCvx_cdc16,SCvx_2018arXiv}, and has been observed to work well in practice. Moreover, the results obtained in~\cite{SzmukReynolds2018} provide preliminary evidence of the real-time capabilities of our methodology.


\newpage


\section{Numerical Results} \label{sec:numerical_results}

In this section we present simulation results that demonstrate the utility of state-triggered and compound state-triggered constraints. Section~\ref{sec4:ex1} presents an example that uses the velocity-triggered angle of attack constraint discussed in~\sref{sec2:stc_ex} to alleviate aerodynamic loading, while~\sref{sec4:ex2} presents an example that uses the collision avoidance compound state-triggered constraint discussed in~\sref{sec2:cmpstc_ex}.

The problem parameters used in the simulations are tabulated in Tables~\ref{tab:params_ex} and 3. All appropriate quantities are expressed in non-dimensionalized mass ($\MU$), length ($\LU$), and time ($\TU$) units. The parameters listed in the tables are drawn from Problem~\ref{prob:ncvx} outlined in~\sref{sec3:ncvx} and from Problem 1, Problem 2, and Algorithm 1 outlined in~\cite{SzmukReynolds2018}. The maximum coast times and the initial positions and velocities are specified individually in~\sref{sec4:ex1} and~\sref{sec4:ex2}. The remaining boundary conditions match those given in the \textit{Boundary Conditions} section of Problem~\ref{prob:ncvx}. Lastly, both simulations assumed the spherical aerodynamics model detailed in Section II.E of~\cite{SzmukReynolds2018}, with the quantity~$\density\Sa\Ca = 0.2$. 

\newcommand{\tabspace}{\hspace{-0.35cm}}
\newcommand{\tabsep}{\null\hspace{0.5cm}\null}
\begin{table}[h]
  \centering
  \caption{Traj. Parameters \hspace{2.7cm} Table 3: Vehicle \& Algorithm Parameters}
  \begin{tabular}{lrllclrll}
  	\hhline{====~====}
    \textbf{Param.}      &   & \tabspace\textbf{Value}                                 & \textbf{Units} &\tabsep&
    \textbf{Param.}      &   & \tabspace\textbf{Value}                                 & \textbf{Units} \\
    \hhline{----~----}
    $\gI$                   &$-$& \tabspace$\ex$                                          & $\LU/\TU^2$ &\tabsep&
	$\rCPB$                 &   & \tabspace$\bvec{0}_{3\times 1}$                         & $\LU$ \\
	$\aoamax$               &   & \tabspace$10.0$                                         & $\dg$ &\tabsep&
    $\rTB$                  &$-$& \tabspace$0.01\cdot\ex$                                 & $\LU$ \\
	$\Vaoa$                 &   & \tabspace$2.5$                                          & $\LU/\TU$ &\tabsep&
    $\mdotalpha$            &   & \tabspace$0.05$                                         & $\TU/\LU$ \\
    $\mig$                  &   & \tabspace$4.0$                                          & $\MU$ &\tabsep&
	$\mdotbeta$             &   & \tabspace$0.02$                                         & $\MU/\TU$ \\
	$\mdry$                 &   & \tabspace$2.0$                                          & $\MU$  &\tabsep&
	$\inertia$              &   & \tabspace$0.01\cdot\diag{\big[0.1\;1\;1\big]}$ & $\MU\cdot\LU^2$ \\
    $\tiltmax$              &   & \tabspace$90.0$                                         & $\dg$ &\tabsep&
    $\KK$                   &   & \tabspace$30$                                           & - \\
	$\omegamax$             &   & \tabspace$90.0$                                         & $\dg/\TU$ &\tabsep&
	$\wvc$                  &   & \tabspace$1\text{e+}3$                                  & - \\
	$\glideslope$           &   & \tabspace$20.0$                                         & $\dg$ &\tabsep&
	$\Wtr$                  &   & \tabspace$\diag{[3\text{e-}6\cdot\bvec{1}^T_{2\times 1}\;1\text{e-}3\cdot\bvec{1}^T_{17\times 1}]}$                                          & - \\
	$\gimbalmax$            &   & \tabspace$20.0$                                         & $\dg$ &\tabsep&
	$\epsvc$                &   & \tabspace$1\text{e-}4$                                  & - \\
    $\Tmin$                 &   & \tabspace$1.0$                                          & $\MU\cdot\LU/\TU^2$ &\tabsep&
    $\epstr$                &   & \tabspace$1\text{e-}1$                                  & - \\
    $\Tmax$                 &   & \tabspace$8.0$                                          & $\MU\cdot\LU/\TU^2$ &\tabsep&
    $\sso$                  &   & \tabspace$10.0$                                         & - \\
	\hhline{====~====}
    \label{tab:params_ex}
  \end{tabular}
\end{table}

\subsection{Example 1: State-Triggered Constraints} \label{sec4:ex1}

This example illustrates the velocity-triggered angle of attack constraint described in~\sref{sec2:stc_ex}. Figures~\ref{fig:traj_aoa_stc_without} and~\ref{fig:traj_aoa_stc_with} show a powered descent guidance trajectory that begins on the free fall trajectory shown in black. The black dot in each of the figures indicates the initial position given by~$\rIin=[14\;16\;0]^T\,\LU$, and is associated with the corresponding initial velocity~$\vIin=[0\;\,-3.57\;\,1.79]^T\;\LU/\TU$. The vehicle is permitted to coast down this path for a maximum duration of~$\tcmax=2.0\;\TU$, at which point the coast phase ends, the engine ignites, and the burn phase commences. The only difference between the scenarios illustrated in Figures~\ref{fig:traj_aoa_stc_without} and~\ref{fig:traj_aoa_stc_with} is that the trajectory of the latter is subjected to the velocity-triggered angle of attack constraint given in~\eqref{eq:stc_ex_aoa_imp}, whereas that of the former is not.

In Figure~\ref{fig:traj_aoa_stc_without}, the vehicle is seen maneuvering away from the free fall trajectory, and uses a bang-bang style thrust profile to get to the landing site located at the origin. The horizontal projections shown at the bottom of Figures~\ref{fig:traj_aoa_stc_without_1} and~\ref{fig:traj_aoa_stc_without_2} clearly show that the vehicle begins its maneuver away from the free fall trajectory as soon as the burn phase begins, without waiting for the vehicle's speed to drop below~$\Vaoa$.

In contrast, the trajectories in Figure~\ref{fig:traj_aoa_stc_with} show the effect of including the velocity-triggered angle of attack constraint. The horizontal projections at the bottom of Figures~\ref{fig:traj_aoa_stc_with_1} and~\ref{fig:traj_aoa_stc_with_2} clearly show that the vehicle travels along the free fall trajectory for a period of time before making a sudden turn towards the landing pad. This behavior is a result of the state-triggered constraint, which requires the vehicle to slow down to a speed below~$\Vaoa$ prior to initiating its large angle of attack maneuver.

\begin{figure}[t!]
	\begin{subfigure}{0.49\textwidth}
    	\centering
		\includegraphics[width=\textwidth]{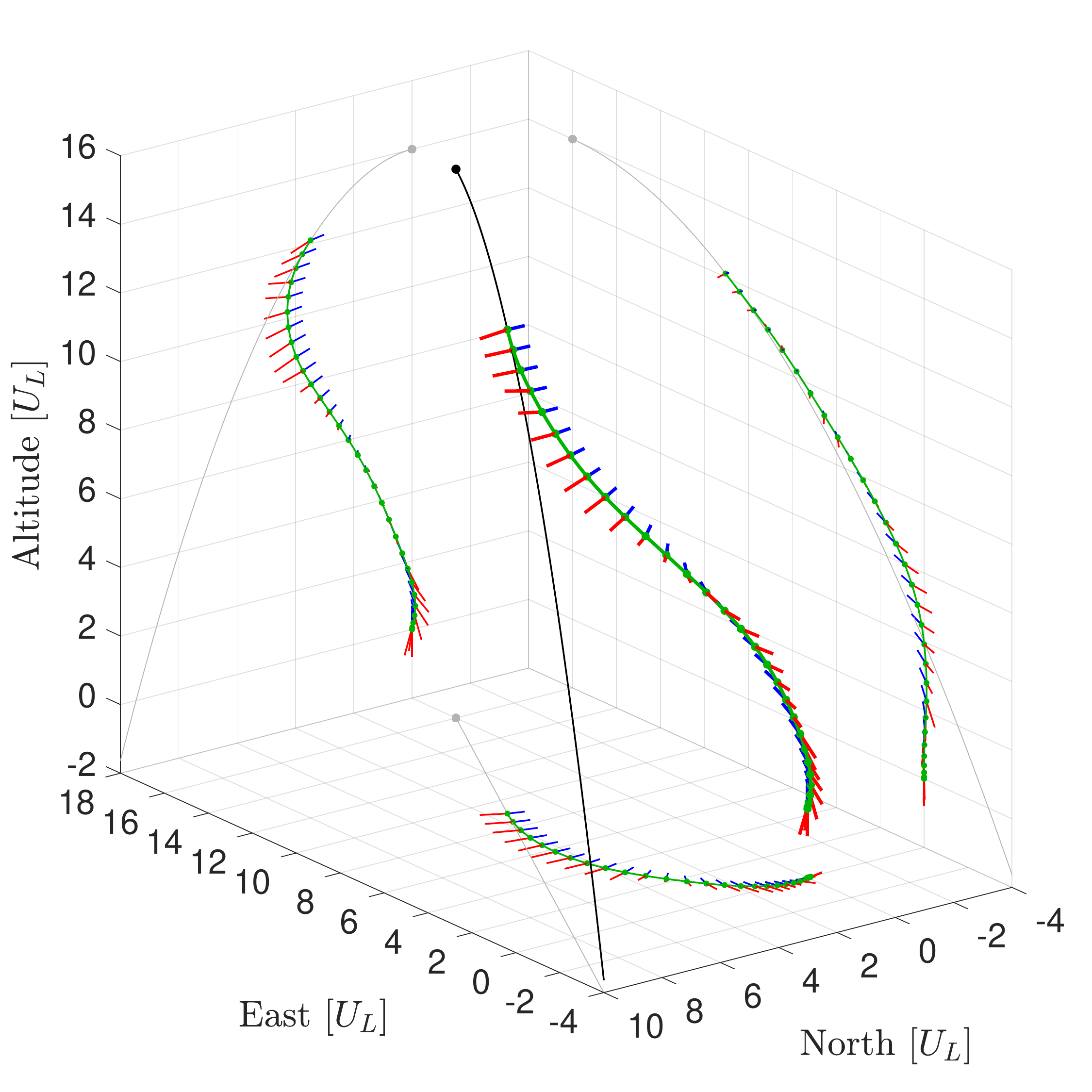}
    	\caption{Up-Range View}
    	\label{fig:traj_aoa_stc_without_1}
   	\end{subfigure}
    \begin{subfigure}{0.49\textwidth}
    	\centering
		\includegraphics[width=\textwidth]{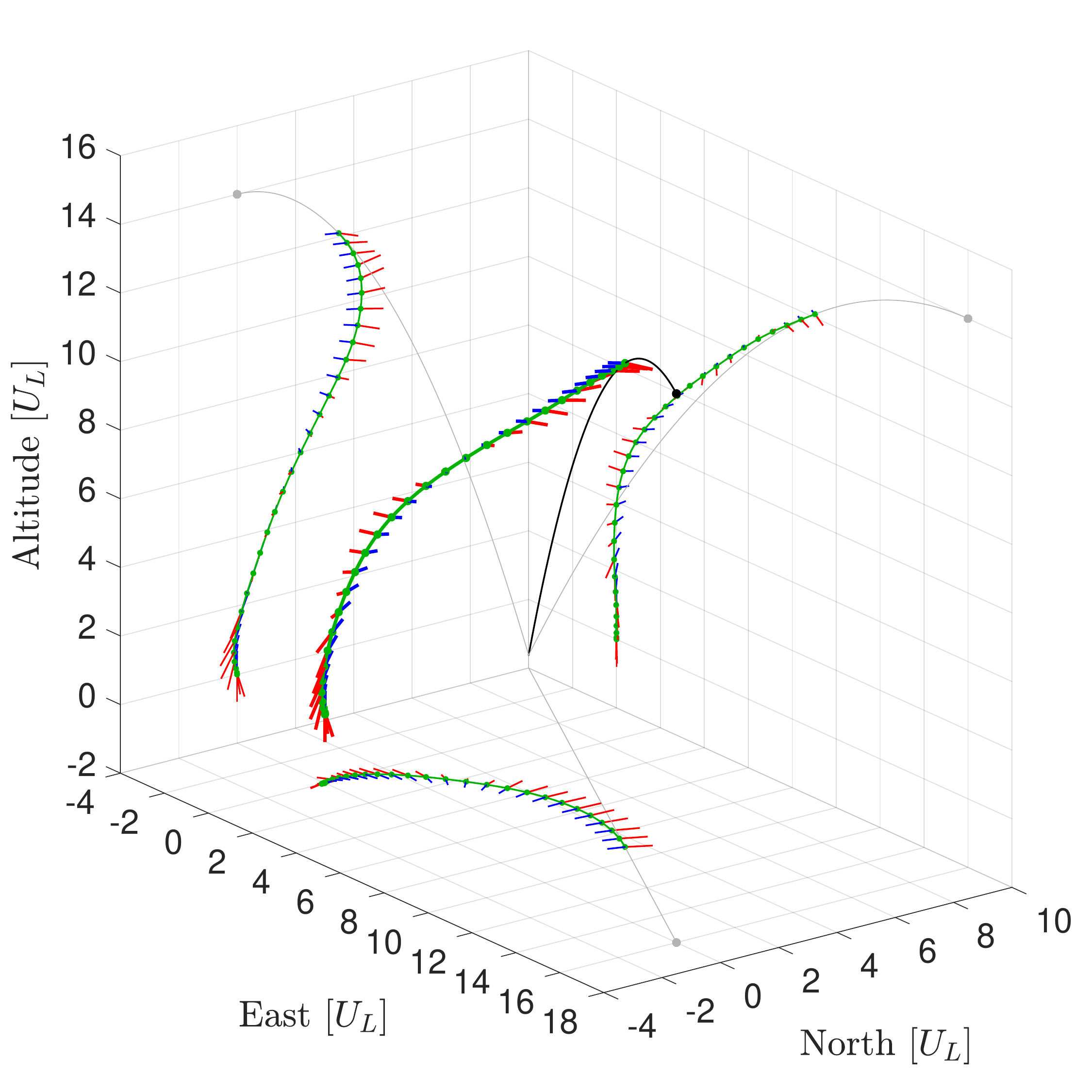}
    	\caption{Down-Range View}
    	\label{fig:traj_aoa_stc_without_2}
   	\end{subfigure}
   	\caption{A powered descent guidance scenario with a free-ignition-time free fall trajectory (shown in black) that \textit{does not} include a state-triggered constraint. The burn phase of the trajectory is shown in green, the thrust vectors in red, and the long axis of the vehicle in blue. Each vector corresponds to a temporal discretization node. Projections of the three-dimensional trajectory are shown on the sides and bottom of the grid box. The two subplots show the same trajectory from different angles.}
   	\label{fig:traj_aoa_stc_without}
\end{figure}

\begin{figure}[t!]
   	\begin{subfigure}{0.49\textwidth}
    	\centering
		\includegraphics[width=\textwidth]{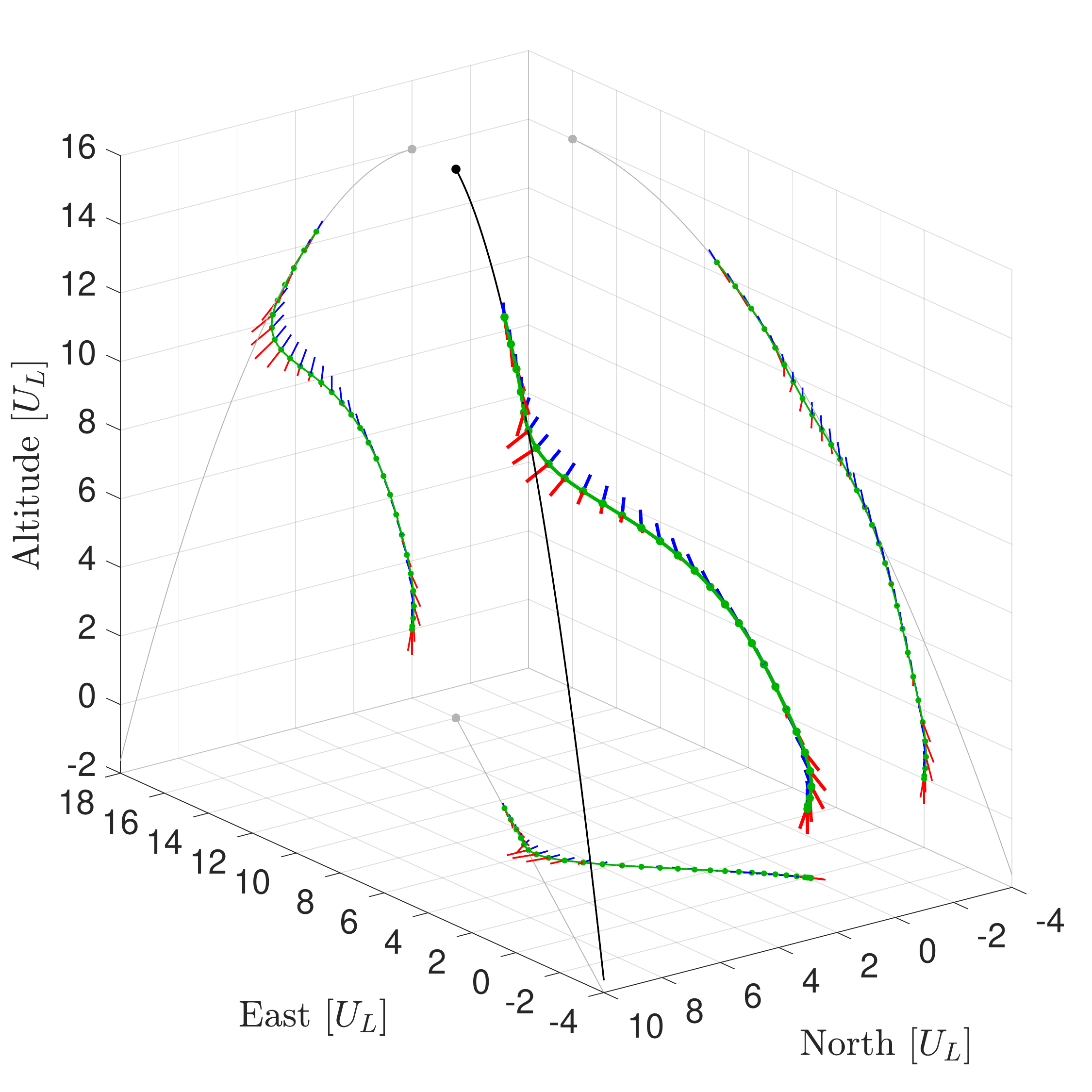}
    	\caption{Up-Range View}
    	\label{fig:traj_aoa_stc_with_1}
   	\end{subfigure}
    \begin{subfigure}{0.49\textwidth}
    	\centering
		\includegraphics[width=\textwidth]{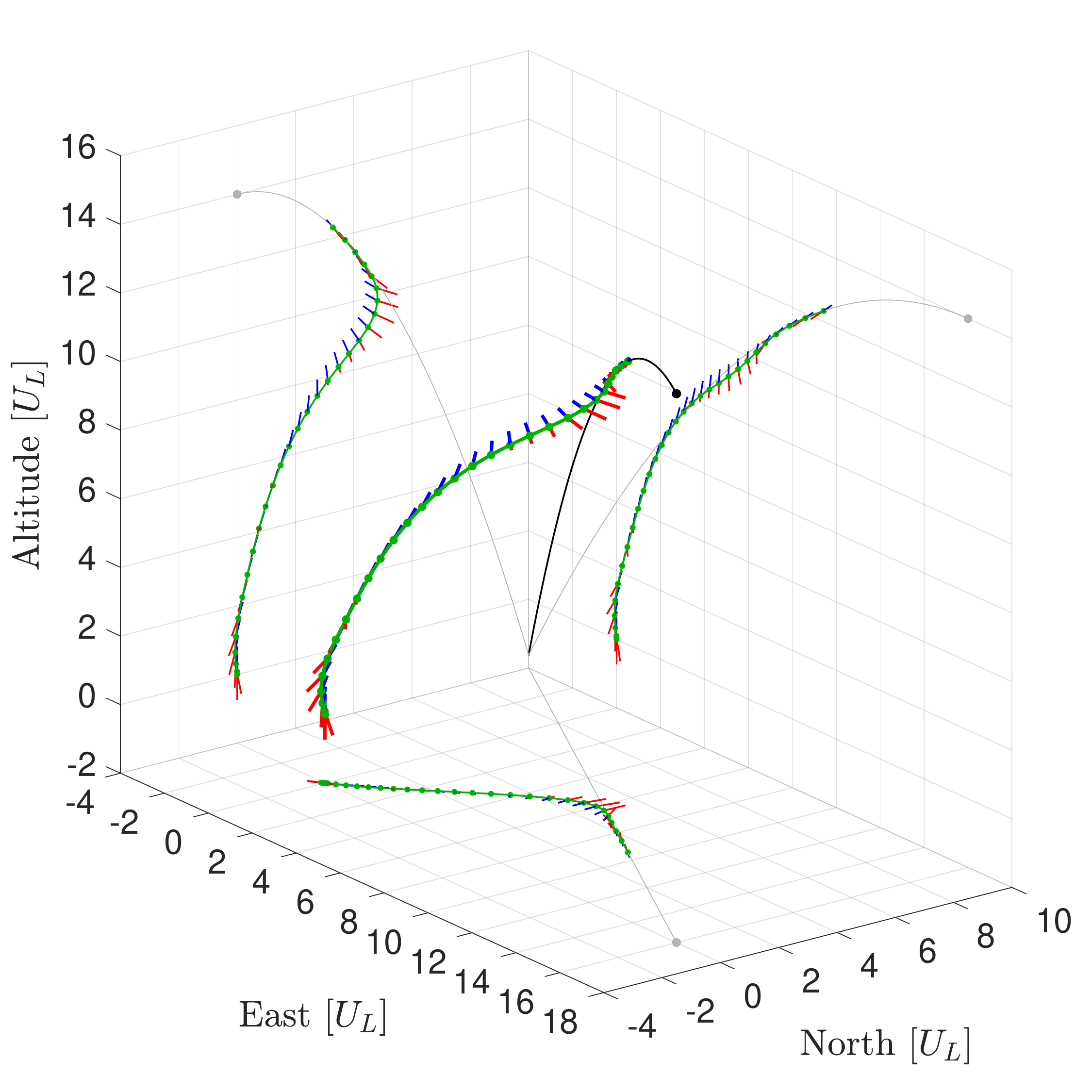}
    	\caption{Down-Range View}
    	\label{fig:traj_aoa_stc_with_2}
   	\end{subfigure}
    \caption{The same scenario as depicted in Figure~\ref{fig:traj_aoa_stc_without}, but with the velocity-triggered angle of attack constraint given in~\eqref{eq:stc_ex_aoa_imp} included.}
   	\label{fig:traj_aoa_stc_with}
   	\vspace{-0.25cm}
\end{figure}

This behavior is shown more explicitly in~Figure~\ref{fig:stc_aoa_stc}, which compares the speed and angle of attack time histories of the two cases. In Figure~\ref{fig:stc_aoa_stc_without} (which corresponds to the trajectories in Figures~\ref{fig:traj_aoa_stc_without_1} and~\ref{fig:traj_aoa_stc_without_2}), the angle of attack clearly violates the prescribed limit of~$\aoamax=10\dg$ when the speed is above the trigger velocity of~$\Vaoa=2.5\;\LU/\TU$. However, in Figure~\ref{fig:stc_aoa_stc_with} (which corresponds to the trajectories in Figures~\ref{fig:traj_aoa_stc_with_1} and~\ref{fig:traj_aoa_stc_with_2}), the behavior is markedly different. Here, the initial angle of attack is less than~$\aoamax$. The vehicle is seen shedding speed between times~$0\;\TU$ and~$1\;\TU$ in an attempt to deactivate the trigger condition. Once it does so, the constraint condition is no longer enforced, and the angle of attack is allowed to increase to close to~$90\dg$. Enforcement of the contrapositive of~\eqref{eq:stc_ex_aoa_orig} can be seen between times~$2\;\TU$ and~$4\;\TU$, where the vehicle is not allowed to break the prescribed speed limit~$\Vaoa$ since the angle of attack is so large. Between times~$4\;\TU$ and~$6\;\TU$, the velocity climbs above~$\Vaoa$ because the angle of attack is less than~$\aoamax$.

This behavior demonstrates that, unlike the multi-phase optimization approach described in Alternative 3 of~\sref{sec2:stc_alt_forms}, the state-triggered constraint is capable of introducing a new constraint phase (i.e. the high speed, low angle of attack phase between times~$4\;\TU$ and~$6\;\TU$) into the solution. This ability allows the vehicle to gain speed during the low angle of attack between times~$4\;\TU$ and~$6\;\TU$, thus attaining a more optimal burn time. Lastly, we note that the case associated with Figures~\ref{fig:traj_aoa_stc_without} and~\ref{fig:stc_aoa_stc_without} is less constrained than the case associated with Figures~\ref{fig:traj_aoa_stc_with} and~\ref{fig:stc_aoa_stc_with}, since the former includes all of the constraints present in the latter except the state-triggered constraint given in~\eqref{eq:stc_ex_aoa_imp}. Thus, as one would expect, the former case has a more optimal (i.e. shorter) burn time when compared to the latter case.

\begin{figure}[t!]
	\begin{subfigure}{0.49\textwidth}
    	\centering
		\includegraphics[width=\textwidth]{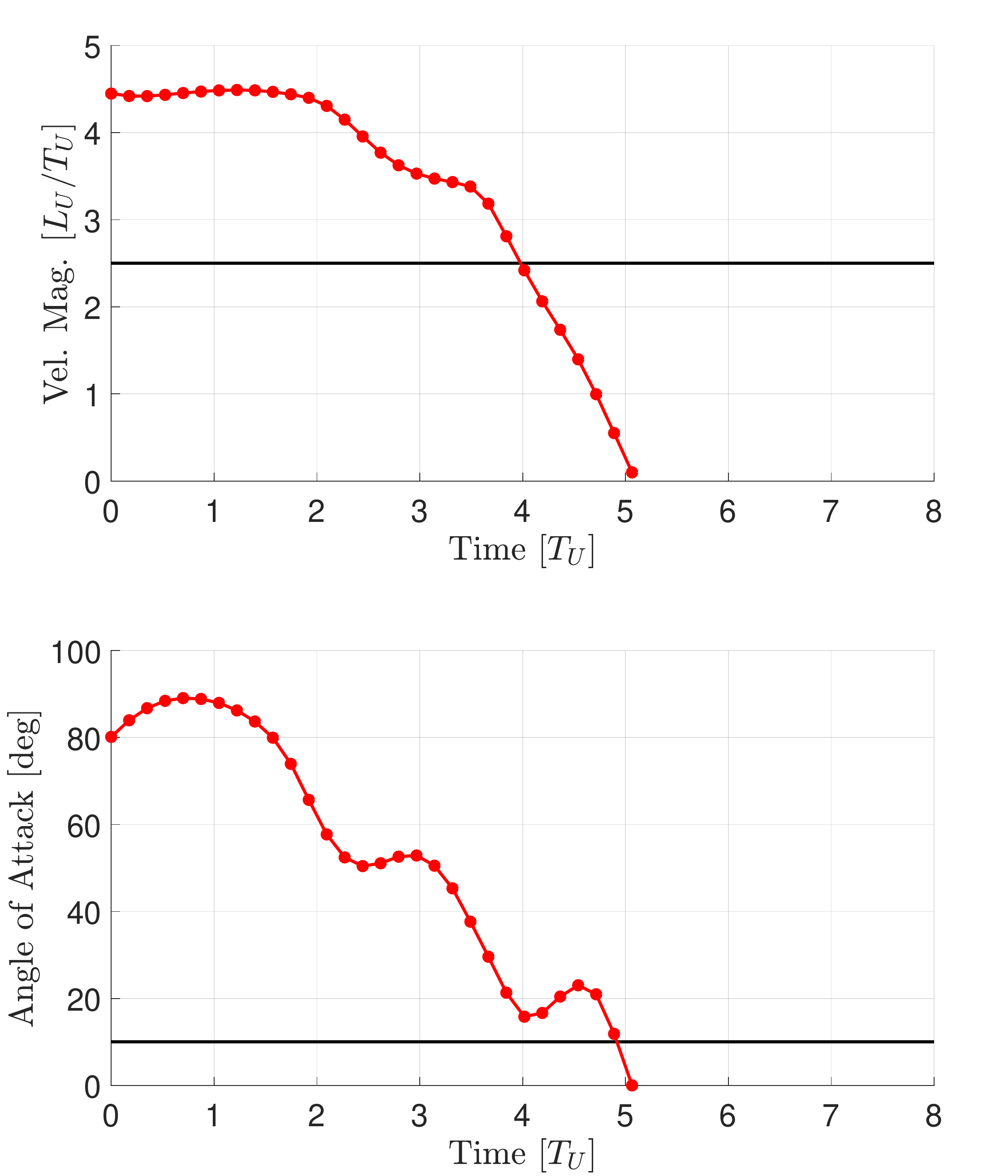}
    	\caption{Without State-Triggered Constraint}
    	\label{fig:stc_aoa_stc_without}
   	\end{subfigure}
    \begin{subfigure}{0.49\textwidth}
    	\centering
		\includegraphics[width=\textwidth]{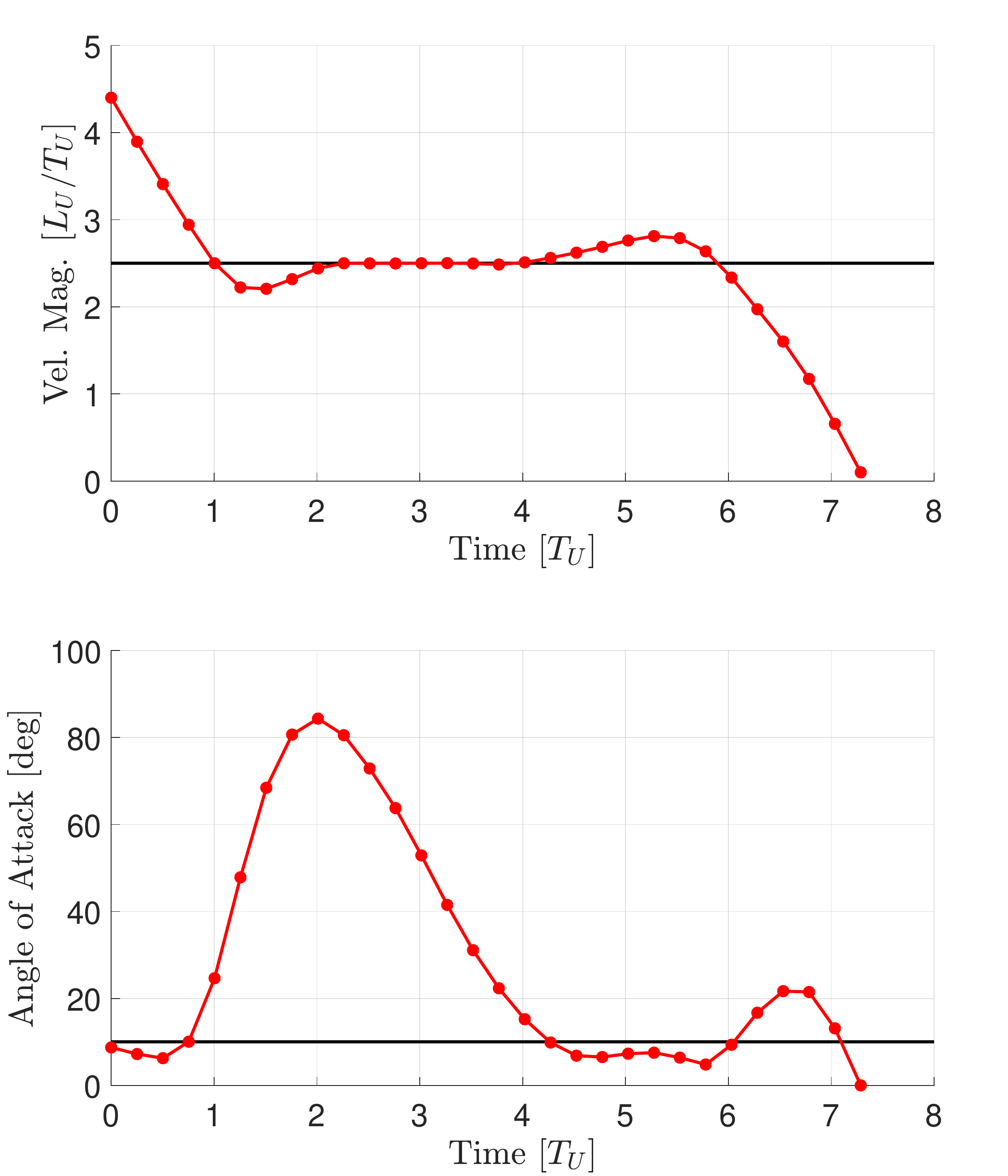}
    	\caption{With State-Triggered Constraint}
    	\label{fig:stc_aoa_stc_with}
   	\end{subfigure}
    \caption{Velocity and angle of attack time histories of the two cases discussed in Example 1. The two figures on the left correspond to the trajectories shown in Figure~\ref{fig:traj_aoa_stc_without}, and the two on the right correspond to the trajectories shown in Figure~\ref{fig:traj_aoa_stc_with}. The black lines in the top two figures represent the trigger velocity~$\Vaoa$, while those in the bottom two figures represent the angle of attack limit~$\aoamax$.}
   	\label{fig:stc_aoa_stc}
   	\vspace{-0.25cm}
\end{figure}

\subsection{Example 2: Compound State-Triggered Constraints} \label{sec4:ex2}

This example illustrates the collision avoidance compound state-triggered constraints discussed in~\sref{sec2:cmpstc_ex}. The scenario used in this example differs from the one used in~\sref{sec4:ex1} in that: (i) the initial position is given by~$\rIin=[6\;10\;0]^T\;\LU$, (ii) the initial velocity is given by~$\vIin=[-1.37\;-3.36\;1.68]^T\;\LU/\TU$, (iii) the velocity-triggered angle of attack constraint given in~\eqref{eq:stc_ex_aoa_imp} \textit{is not} included, and (iv) a cubic keep-out volume is included in the vicinity of the landing site. For presentation purposes, the keep-out volume is a cube with four walls, instead of the two-walled volume presented in~\sref{sec2:cmpstc_ex}. This keep-out volume requires that the compound state-triggered constraint in~\eqref{eq:cmpstc_ex_and} be augmented with two more trigger conditions. This modification is rather easy to make, and is assumed in the following discussion.

Figures~\ref{fig:cmpstc_messa_without} and~\ref{fig:cmpstc_messa_with} show the trajectories for the scenario described above. Figures~\ref{fig:cmpstc_messa_without_1} and~\ref{fig:cmpstc_messa_without_2} show two views of the same trajectory, which was generated \textit{without} the inclusion of~\eqref{eq:cmpstc_ex_and}. This trajectory is seen to violate the keep-out volume. This violation is most evident in the far-right projection in Figures~\ref{fig:cmpstc_messa_without_1} and~\ref{fig:cmpstc_messa_without_2}.

In contrast, Figures~\ref{fig:cmpstc_messa_with_1} and~\ref{fig:cmpstc_messa_with_2}, show the trajectory computed when the collision avoidance compound state-triggered constraint given \textit{is} included. In this case, the vehicle is able to avoid violations of the keep-out volume at each of the temporal nodes, as seen in the far-right projections shown in Figures~\ref{fig:cmpstc_messa_with_1} and~\ref{fig:cmpstc_messa_with_2}. Note that the continuous-time trajectory does violate the keep-out volume between the~$\ith{14}{th}$ and~$\ith{15}{th}$ temporal nodes. This ``clipping'' behavior is a common artifact of discretization, and can be addressed in a variety of ways. However, we deem such considerations to be beyond the scope of this paper.

This example demonstrates that our compound state-triggered constraints can be incorporated into a powered-descent guidance problem with relative ease. These constraints allow the designer to construct complicated keep-out geometries using relatively simple half-space constraints. In setting up this example, we found that the proposed compound state-triggered constraint formulation worked well for a variety of initial conditions.

\begin{figure}[t!]
	\begin{subfigure}{0.49\textwidth}
    	\centering
		\includegraphics[width=\textwidth]{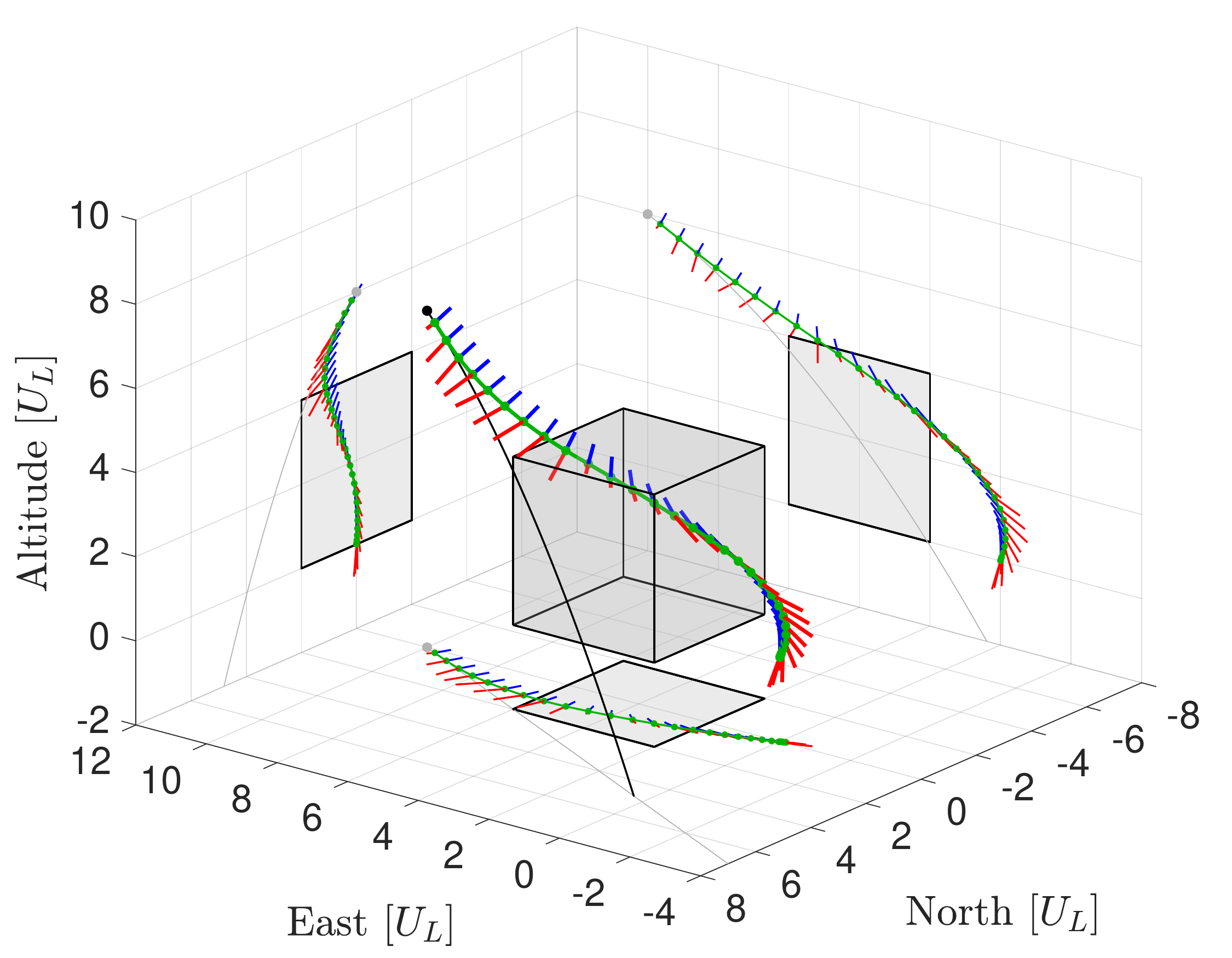}
    	\caption{Up-Range View}
    	\label{fig:cmpstc_messa_without_1}
   	\end{subfigure}
    \begin{subfigure}{0.49\textwidth}
    	\centering
		\includegraphics[width=\textwidth]{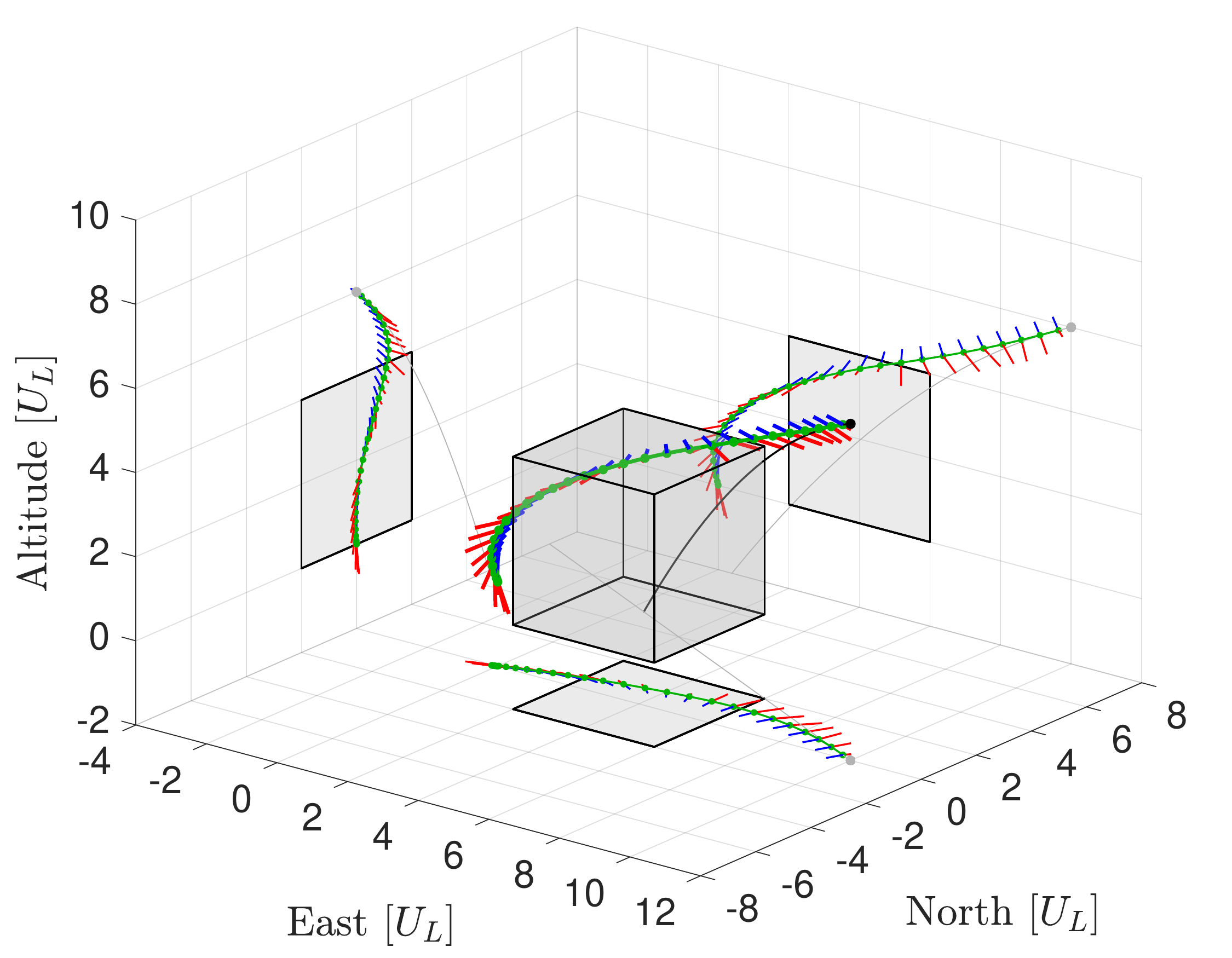}
    	\caption{Down-Range View}
    	\label{fig:cmpstc_messa_without_2}
   	\end{subfigure}
    \caption{A powered descent guidance scenario with a cubic keep-out volume located near the landing site. In this case, the trajectory was computed \textit{without} the inclusion of a collision avoidance constraint. The two figures present the same trajectory shown from different angles. The trajectory is seen to violate the keep-out volume in the far-right projection shown in both views.}
   	\label{fig:cmpstc_messa_without}
\end{figure}

\begin{figure}[t!]
	\begin{subfigure}{0.49\textwidth}
    	\centering
		\includegraphics[width=\textwidth]{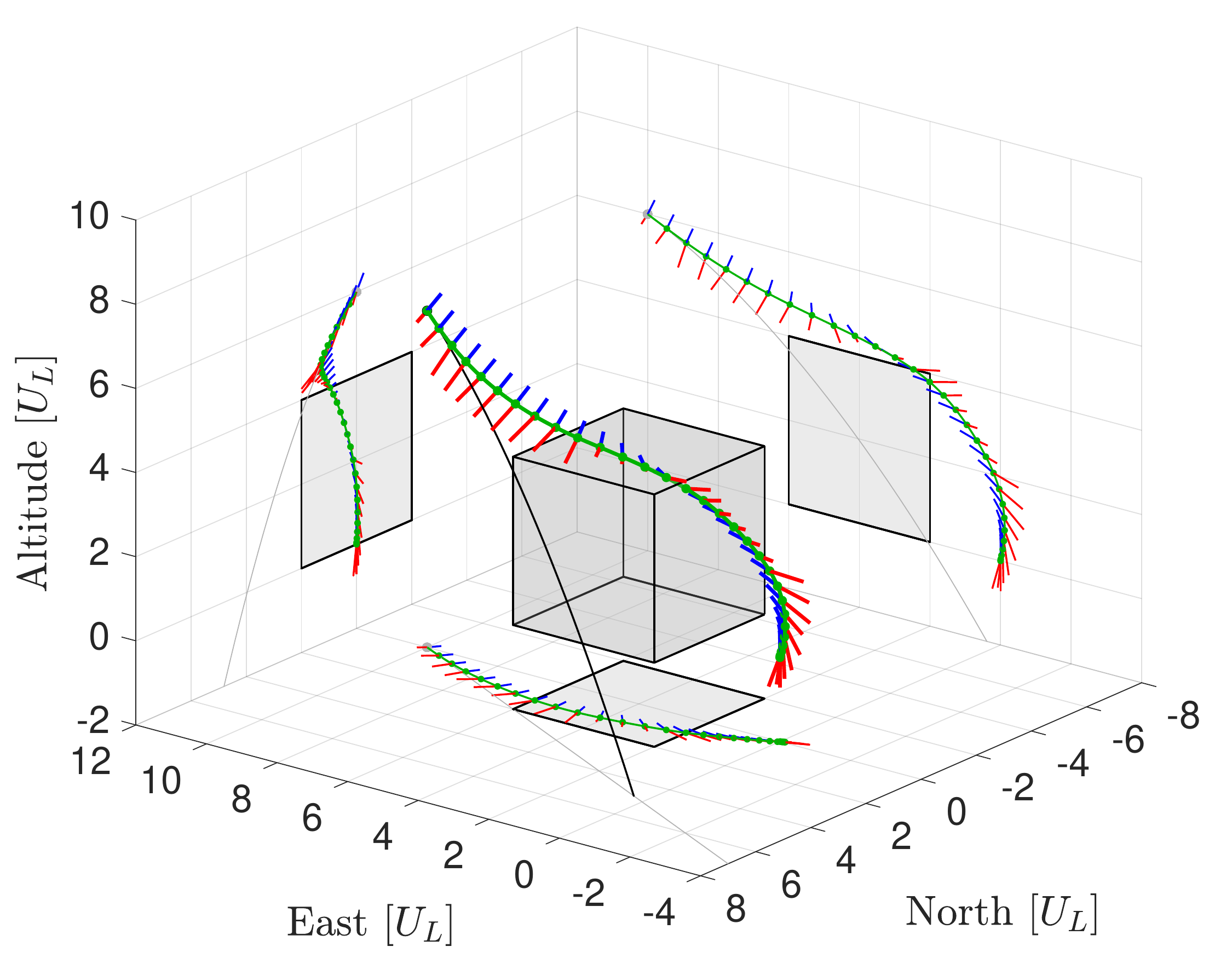}
    	\caption{Up-Range View}
    	\label{fig:cmpstc_messa_with_1}
   	\end{subfigure}
    \begin{subfigure}{0.49\textwidth}
    	\centering
		\includegraphics[width=\textwidth]{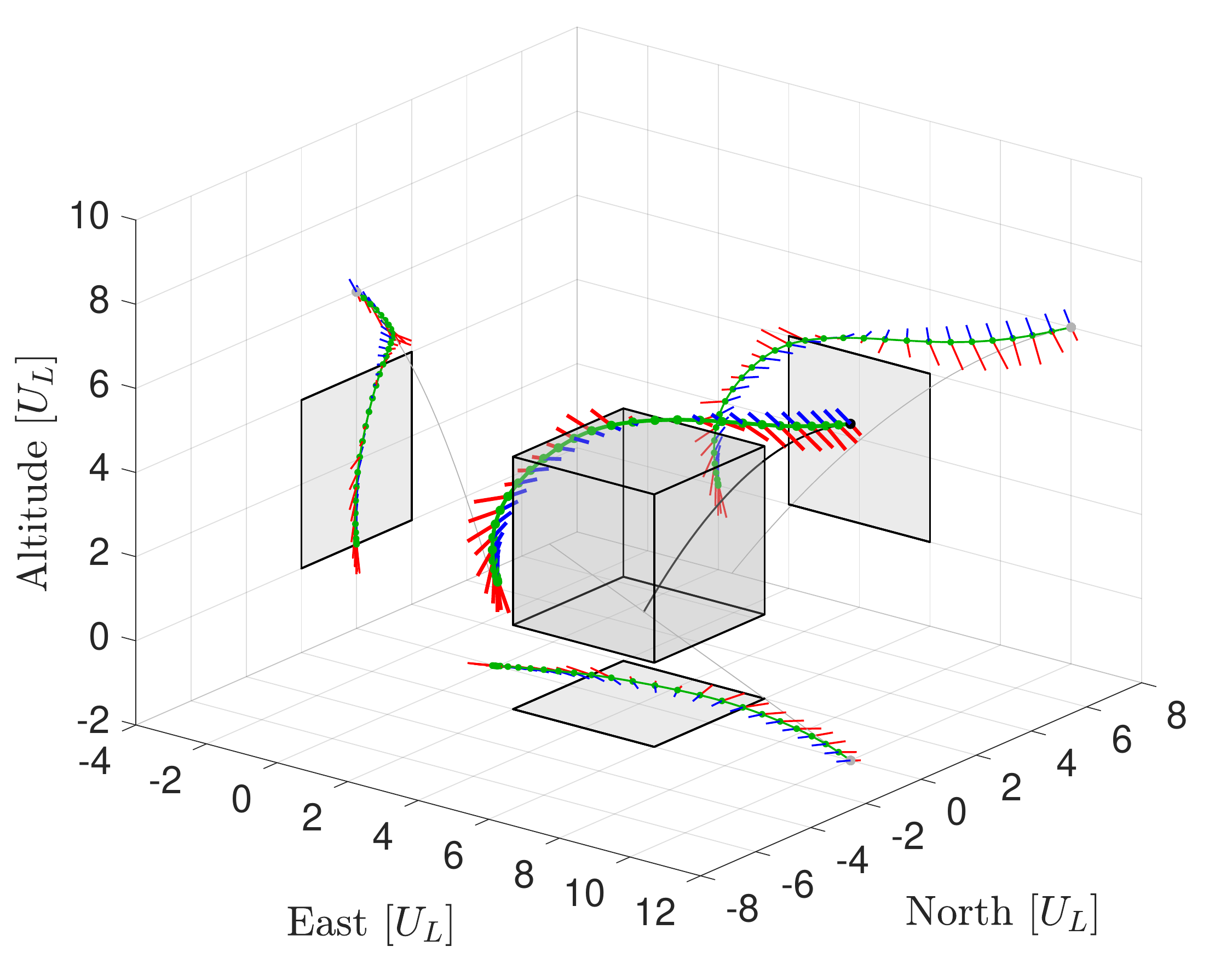}
    	\caption{Down-Range View}
    	\label{fig:cmpstc_messa_with_2}
   	\end{subfigure}
    \caption{The scenario corresponding to the one depicted in Figure~\ref{fig:cmpstc_messa_without} \textit{with} the inclusion of a collision avoidance compound state-triggered constraint given in~\eqref{eq:cmpstc_ex_and}. Here, the trajectory steers clear of the keep-out volume, as seen in the far-right projection seen in both views.}
   	\label{fig:cmpstc_messa_with}
  	\vspace{-0.25cm}
\end{figure}


\section{Conclusion} \label{sec:conclusion}

This paper presents a powered descent guidance problem formulation that contains state-triggered and compound state-triggered constraints. The primary contribution of this paper is a novel continuous formulation for \textit{compound} state-triggered constraints. This formulation is a generalization of the state-triggered constraint formulation presented in~\cite{SzmukReynolds2018}, and enables the use of trigger and constraint conditions composed using Boolean \textit{and} and \textit{or} operations. Both state-triggered constraint formulations are demonstrated in example simulations using the successive convexification framework. The examples showcase a velocity-triggered angle of attack constraint that is used to alleviate aerodynamic loads, and a collision avoidance constraint used to avoid large geological formations near the landing site. Notably, the velocity-triggered angle of attack constraint demonstrates the ability of state-triggered constraints to introduce and rearrange constraint phases in the solution, without resorting to combinatorial techniques. While the results presented here and in~\cite{SzmukReynolds2018} show promise for powered descent guidance applications, future work will focus on improving the convergence behavior of problems including state-triggered and compound state-triggered constraints.


\section{Acknowledgements}

This research has been supported by NASA grant NNX17AH02A and was partially carried out at the Johnson Space Center. Government sponsorship acknowledged. Support for studying successive convexification was provided by the Office of Naval Research grants N00014-16-1-2877 and N00014-16-3144.


\bibliography{bibliography}

\end{document}